\documentclass[10pt]{amsart}

\usepackage{amsmath, amssymb, euscript}
\usepackage{mathscinet}
\usepackage[all]{xy}

\theoremstyle{plain}
\newtheorem{thm}{Theorem}
\newtheorem{prop}[thm]{Proposition}
\newtheorem{lem}[thm]{Lemma}
\newtheorem{cor}[thm]{Corollary}
\newtheorem{ithm}{Theorem}

\theoremstyle{remark}

\newtheorem{exc}{Example}

\newcommand{\Acat}[1]{\mathsf{A}_{#1}}
\newcommand{\acat}[1]{\mathsf{a}_{#1}}
\newcommand{\cacat}[1]{\mathsf{\hat{a}}_{#1}}
\newcommand{\lcat}{\mathsf{l}}
\newcommand{\sets}{\mathsf{Sets}}

\newcommand{\cat}[1]{\mathsf{#1}}
\newcommand{\ccat}{\mathsf{c}}
\newcommand{\morc}{\mathsf{Mor \; c}}
\newcommand{\morcat}[1]{\mathsf{Mor \; {#1}}}
\newcommand{\mcat}[1]{\mathsf{Mod}(#1)}
\newcommand{\smcat}[1]{\mathsf{mod}(#1)}
\newcommand{\presh}{\mathsf{PreSh}}
\newcommand{\shc}[1]{\mathsf{Sh}(X,#1)}

\newcommand{\qcoh}[1]{\mathsf{QCoh}(#1)}
\newcommand{\compl}{\mathsf{Compl}}
\newcommand{\oc}[1]{\mathsf{#1}}

\newcommand{\sh}[1]{\mathcal{#1}}
\newcommand{\sa}{\mathcal{A}}
\newcommand{\sfam}{\mathcal{F}}
\newcommand{\mfam}[1]{\mathcal{#1}}
\newcommand{\defs}{\mathsf{Def}_{\sfam}}
\newcommand{\defm}[1]{\mathsf{Def}_{#1}}
\newcommand{\df}{\mathsf{D}}
\newcommand{\hull}{\mathsf{H}}
\newcommand{\infdef}{\mathit{T}^1}
\newcommand{\obstrdef}{\mathit{T}^2}
\newcommand{\qc}[1]{\mathsf{qc}_{#1}}

\newcommand{\shder}{\mathcal{D}er}
\newcommand{\shhom}{\mathcal{H}om}
\newcommand{\shend}{\mathcal{E}nd}
\newcommand{\shext}{\mathcal{E}xt}

\newcommand{\ch}{\mathsf{H}}
\newcommand{\cpd}{*}
\newcommand{\hc}{\mathsf{HC}}
\newcommand{\hh}{\mathsf{HH}}
\newcommand{\cc}{\mathsf{C}}
\newcommand{\dc}{\mathsf{D}}
\newcommand{\dI}[1]{d^{#1}_{I}}
\newcommand{\dII}[1]{d^{#1}_{II}}

\newcommand{\Gm}{\Gamma}

\newcommand{\reg}[1]{\mathcal{O}_{#1}}
\newcommand{\diff}[1]{\mathrm{Diff}(\sh{#1})}
\newcommand{\g}{\mathbf{g}}
\newcommand{\ue}{\textup{U}}

\newcommand{\pp}{\mathbf{P}}
\newcommand{\rst}[1]{|_{#1}}

\DeclareMathOperator{\car}{char}
\DeclareMathOperator{\coker}{coker}
\DeclareMathOperator{\der}{Der}

\DeclareMathOperator{\ext}{Ext}
\DeclareMathOperator{\gr}{gr}
\DeclareMathOperator{\hmm}{Hom}
\DeclareMathOperator{\id}{id}

\DeclareMathOperator{\mor}{Mor}
\DeclareMathOperator*{\osum}{\oplus}

\title[Noncommutative deformations]{Computing noncommutative deformations \\
of presheaves and sheaves of modules}
\author{Eivind Eriksen}
\address{Oslo University College \\ Postboks 4 St. Olavs Plass \\
0130 Oslo, Norway} \email{eeriksen@hio.no}

\thanks{This research has been supported by a postdoctoral grant awarded
by the Norwegian Research Council, project no. 157740/432, and a grant
awarded by the Mittag-Leffler Institute.}

\begin{document}

\begin{abstract}
We describe a noncommutative deformation theory for presheaves and
sheaves of modules that generalizes the commutative deformation
theory of these global algebraic structures, and the noncommutative
deformation theory of modules over algebras due to Laudal.

In the first part of the paper, we describe a noncommutative
deformation functor for presheaves of modules on a small category, and
an obstruction theory for this functor in terms of global Hochschild
cohomology. An important feature of this obstruction theory is that it
can be computed in concrete terms in many interesting cases.

In the last part of the paper, we describe noncommutative
deformation functors for sheaves and quasi-coherent sheaves of
modules on a ringed space $(X,\sa)$. We show that for any good
$\sa$-affine open cover $\oc U$ of $X$, the forgetful functor $\qcoh
\sa \to \presh(\oc U, \sa)$ induces an isomorphism of noncommutative
deformation functors.

\emph{Applications.} We consider noncommutative deformations of
quasi-coherent $\sa$-modules on $X$ when $(X, \sa) = (X, \reg X)$ is
a scheme or $(X, \sa) = (X, \sh D)$ is a D-scheme in the sense of
Beilinson and Bernstein. In these cases, we may use any open affine
cover of $X$ closed under finite intersections to compute
noncommutative deformations in concrete terms using presheaf
methods. We compute the noncommutative deformations of the left $\sh
D_X$-module $\reg X$ when $X$ is an elliptic curve as an example.
\end{abstract}

\maketitle

\section*{Introduction}

Deformation theory was formalized by Grothendieck in the language of
schemes in the 1950s, and is described in the series of Bourbaki
seminar expositions \emph{Fondements de la g\'eom\'etrie
alg\'ebrique} \cite{grot62}; see in particular Grothendieck
\cite{grot95i}, \cite{grot95ii}. The general philosophy is best
described by the following quotation:
\begin{quote}
La m\'ethode g\'en\'erale consiste toujours \`a faire des
constructions formelles, ce qui consiste essentiellement \`a faire
de la g\'eom\'etrie alg\'ebrique sur un anneau artinien, et \`a en
tirer des conclusions de nature "alg\'ebrique" en utilisant les
trois th\'eor\`emes fondamentaux (Grothendieck \cite{grot95i}, p.
11)
\end{quote}
We shall follow Grothendieck's philosophy closely, and we are therefore
led to the study of functors of (noncommutative) Artin rings.

Let $k$ be an algebraically closed field, and let $\lcat$ denote the
category of local Artinian commutative $k$-algebras with residue
field $k$, with local homomorphisms. A \emph{functor of Artin rings}
is a covariant functor $\df: \lcat \to \sets$ such that $\df(k)$
only contains one element. In Schlessinger \cite{schl68}, criteria
for functors of Artin rings to have a pro-representable hull,
respectively for functors of Artin rings to be pro-representable,
were given.

Let $\cat A$ be an Abelian $k$-category, and let $X$ be an object of
$\cat A$. The flat deformation functor $\defm X: \lcat \to \sets$ of
$X$ in $\cat A$ is a functor of Artin rings. In many cases, it has a
pro-representing hull $\hull(\defm X)$, see Schlessinger \cite{schl68},
Laudal \cite{laud79}, and there are constructive methods for finding
$\hull(\defm X)$, see Laudal \cite{laud79}, \cite{laud86}. In fact, if
there exists an obstruction theory for $\defm X$ with finite
dimensional cohomology $\ch^n(X)$ for $n=1,2$, then $\defm X$ has a
pro-representing hull, given algorithmically in terms of the vector
spaces $\ch^n(X)$ for $n = 1,2$ and certain generalized symmetric
Massey products on them.

When $\cat A = \mcat A$, the category of left modules over an
associative $k$-algebra $A$, Laudal introduced a generalization of the
deformation functor $\defm M: \lcat \to \sets$ for a left $A$-module
$M$ in Laudal \cite{laud02}. He considered the category $\acat p$ of
$p$-pointed Artinian rings for any integer $p \ge 1$, and constructed a
noncommutative deformation functor $\defm{\mfam M}: \acat p \to \sets$
for any finite family $\mfam M = \{ M_1, \dots, M_p \}$ of left
$A$-modules. This deformation functor has an obstruction theory with
cohomology $( \ext^n_A(M_j,M_i) )$, and a pro-representing hull
$\hull(\defm{\mfam M})$, given algorithmically in terms of the vector
spaces $( \ext^n_A(M_j,M_i) )$ and certain generalized Massey products
on them.

The objects in the category $\acat p$ are Artinian rings $R$, together
with ring homomomorphisms $k^p \to R \to k^p$ such that the composition
is the identity, and such that $R$ is $I$-adic complete for $I = \ker(R
\to k^p)$. The morphisms are the natural commutative diagrams. In
section \ref{s:functors}, we give a systematic introduction to functors
$\df: \acat p \to \sets$ of noncommutative Artian rings, following
Laudal \cite{laud02}.

The idea is that there is a noncommutative deformation functor $\defm
X: \acat p \to \sets$ for any finite family $X = \{ X_1, \dots, X_p \}$
of algebraic or algebro-geometric objects. The restrictions of $\defm
X$ along the $p$ natural full embeddings of categories $\acat 1
\subseteq \acat p$ are the noncommutative deformation functors
$\defm{X_i}: \acat 1 \to \sets$, and the restriction of $\defm{X_i}$ to
$\lcat \subseteq \acat 1$ is the commutative deformation functor
$\defm{X_i}^c: \lcat \to \sets$. When these deformation functors have
pro-representing hulls, we show that
    \[ \hull(\defm X)^\text{comm} \cong \osum_{1 \le i \le p}
    \hull(\defm{X_i})^\text{comm} \cong \osum_{1 \le i \le p}
    \hull(\defm{X_i}^c) \]
We remark that the hull $\hull(\defm X)$ \emph{is not} isomorphic to
$\oplus_i \hull(\defm{X_i})$ in general. We also remark that the
noncommutative deformation functor $\defm X: \acat p \to \sets$ of a
family $X = \{ X_1, \dots, X_p \}$ \emph{is not} the same as the
noncommutative deformation functor $\defm X: \acat 1 \to \sets$ of the
direct sum $X = X_1 \oplus \dots \oplus X_p$.

Let $\presh(\ccat, \sa)$ be the category of presheaves of left
$\sa$-modules on $\ccat$, where $\ccat$ is a small category and
$\sa$ is a $k$-algebra of presheaves on $\ccat$. We consider an
Abelian $k$-category $\cat A$ such that $\cat A \subseteq
\presh(\ccat, \sa)$ is a full subcategory, and construct a
noncommutative deformation functor $\defs^{\cat A}: \acat p \to
\sets$ for any finite family $\sfam = \{ \sh F_1, \dots, \sh F_p \}$
of objects in $\cat A$ in section \ref{s:defm}. We remark that we
use a notion of matric freeness, introduced by Laudal, to replace
flatness in the definition of the deformation functor. It is not
clear if these notions are equivalent if $p \ge 2$.

In section \ref{s:defm-presh} - \ref{s:obstr}, we consider
deformations in the category $\cat A = \presh(\ccat, \sa)$ of
presheaves. We describe the noncommutative deformation functor
$\defs: \acat p \to \sets$ of any finite family $\sfam = \{ \sh F_1,
\dots, \sh F_p \}$ of presheaves of left $\sa$-modules on $\ccat$ in
concrete terms in section \ref{s:defm-presh}, and use this to
develop an obstruction theory for $\defs$ with cohomology $(
\hh^n(\ccat, \sa, \shhom_k(\sh F_j, \sh F_i)) )$ in section
\ref{s:obstr}. The global Hochschild cohomology $\hh^n(\ccat, \sa,
\shhom_k(\sh F_j, \sh F_i))$ of $\sa$ with values in $\shhom_k(\sh
F_j, \sh F_i)$ on $\ccat$ is described in detail in section
\ref{s:coh}.

\begin{ithm}
Let $\ccat$ be a small category, let $\sa$ be a presheaf of
$k$-algebras on $\ccat$, and let $\sfam = \{ \sh F_1, \dots, \sh F_p
\}$ be a finite family of presheaves of left $\sa$-modules on $\ccat$.
If $\dim_k \hh^n(\ccat, \sa, \shhom_k(\sh F_j, \sh F_i)) < \infty
\text{ for } 1 \le i,j \le p, \; n = 1,2$, then the noncommutative
deformation functor $\defs: \acat p \to \sets$ of $\sfam$ in
$\presh(\ccat, \sa)$ has a pro-representing hull $\hull(\defs)$,
completely determined by the $k$-linear spaces $\hh^n(\ccat, \sa,
\shhom_k(\sh F_j, \sh F_i))$ for $1 \le i,j \le p$, $n = 1,2$, together
with some generalized Massey products on them.
\end{ithm}

In section \ref{s:defmsh}, we consider deformations in the category
$\cat A = \shc \sa$ of sheaves of left $\sa$-modules on $X$, where
$(X, \sa)$ is a ringed space over $k$. We may consider $\shc \sa$ as
a full subcategory of $\presh(\ccat(X),\sa)$, where $\ccat(X)$ is
the category with open subsets $U \subseteq X$ as objects, and
opposite inclusion $U \supseteq V$ as morphisms. Moreover, $\shc
\sa$ is an Abelian $k$-category. We show that for any finite family
$\sfam$ in $\shc \sa$, the natural forgetful functor $\shc \sa \to
\presh(\ccat(X), \sa)$ induces an isomorphism $\defs^{sh} \to \defs$
of deformation functors. However, this result is not very useful for
computational purposes, since the category $\ccat(X)$ consisting of
all open sets in $X$ is usually too big to allow for effective
computations.

In section \ref{s:qcoh}, we consider deformations in the category
$\cat A = \qcoh \sa$ of quasi-coherent sheaves of left $\sa$-modules
on $X$, where $(X,\sa)$ is a ringed space over $k$. If $X$ has an
$\sa$-affine open cover $\oc U$, then $\qcoh \sa$ is an Abelian
$k$-category and a full subcategory of $\presh(\ccat(X),\sa)$. We
show that if $\oc U$ is a good $\sa$-affine open cover of $X$, then
the natural forgetful functor $\qcoh \sa \to \presh(\oc U, \sa)$
induces an isomorphism $\defs^{qc} \to \defs^{\oc U}$ of
noncommutative deformation functors for any finite family $\sfam$ in
$\qcoh \sa$.

\begin{ithm}
Let $(X,\sa)$ be a ringed space over $k$, let $\oc U$ be a good
$\sa$-affine open cover of $X$, and let $\sfam = \{ \sh F_1, \dots, \sh
F_p \}$ be a finite family of quasi-coherent left $\sa$-modules on $X$.
If $\dim_k \hh^n(\oc U, \sa, \shhom_k(\sh F_j, \sh F_i)) < \infty$ for
$1 \le i,j \le p, \; n = 1,2$, then the noncommutative deformation
functor $\defs^{qc}: \acat p \to \sets$ of $\sfam$ in $\qcoh \sa$ has a
pro-representing hull $\hull(\defs^{qc})$, completely determined by the
$k$-linear spaces $\hh^n(\oc U, \sa, \shhom_k(\sh F_j, \sh F_i))$ for
$1 \le i,j \le p$, $n = 1,2$, together with some generalized Massey
products on them.
\end{ithm}

We give examples of ringed spaces $(X, \sa)$ that admit good
$\sa$-affine open covers in section \ref{s:imp-ex}. The main
commutative examples are schemes $(X,\reg X)$ over $k$. The main
noncommutative examples are D-schemes $(X, \sh D)$ over $k$ in the
sense of Beilinson, Bernstein \cite{bei-ber93}. In particular,
important examples of D-schemes over an algebraically closed field
$k$ of characteristic $0$ include $(X,\sh D_X)$, where $X$ is a
locally Noetherian scheme over $k$ and $\sh D_X$ is the sheaf of
$k$-linear differential operators on $X$, and $(X,\ue(\g))$, where
$X$ is a separated scheme of finite type over $k$ and $\ue(\g)$ is
the universal enveloping D-algebra of a Lie algebroid $\g$ on $X/k$.

Let us consider one of the Abelian categories $\cat A =
\presh(\ccat, \sa)$, $\cat A = \shc \sa$, and if $X$ has a good
$\sa$-affine open cover, $\cat A = \qcoh \sa$. We expect that the
noncommutative deformation functor $\defs^{\cat A}: \acat p \to
\sets$ of $\sfam$ in $\cat A$ is controlled by $(\ext^n_{\cat A}(\sh
F_j, \sh F_i))$ in all these cases. In fact, we show that
$t(\defs^{\cat A})_{ij} \cong \ext^1_{\cat A}(\sh F_j, \sh F_i)$, at
least if $\dim_k \hh^n(\oc U, \sa, \shhom_k(\sh F_j, \sh F_i)) <
\infty \text{ for } 1 \le i,j \le p, \; n = 1,2$, see corollary
\ref{c:defm-iso}, the remark in section \ref{s:defmsh}, and the
remark following theorem \ref{t:defmqc-hull}. It might be possible
to develop an obstruction theory for $\defs^{\cat A}$ with
cohomology $(\ext^n_{\cat A}(\sh F_j, \sh F_i))$ in all these cases.

However, notice that when $\cat A = \shc \sa$ or $\cat A = \qcoh
\sa$ for a noncommutative sheaf of algebras $\sa$, it is often very
hard to compute $\ext^n_{\cat A}(\sh F_j, \sh F_i)$. For instance,
localization of injectives can behave badly, even when $(X, \sa)$ is
a D-scheme. On the other hand, the global Hochschild cohomology
groups $\hh^n(\oc U, \sa, \shhom_k(\sh F_j, \sh F_i))$ can be
computed in concrete terms in many cases of interest.

In section \ref{s:example}, we give an example of this. Let $X$ be any
elliptic curve over an algebraically closed field $k$ of characteristic
$0$, and consider $\reg X$ as a left $\sh D_X$-module on $X$. We show
that $\hh^0(\oc U, \sh D_X, \shend_k(\reg X)) \cong k$, $\hh^1(\oc U,
\sh D_X, \shend_k(\reg X)) \cong k^2$, and $\hh^2(\oc U, \sh D_X,
\shend_k(\reg X)) \cong k$ for an open affine cover $\oc U$ of $X$
closed under intersections. Using these results and the obstruction
calculus, we compute the pro-representing hull $\hull(\defm{\reg X})
\cong k{\ll}t_1,t_2{\gg}/(t_1 t_2 - t_2 t_1)$ of the noncommutative
deformation functor $\defm{\reg X}$ and we also compute its versal
family. In this example, it seems hard to compute $\ext^n_{\sh
D_X}(\reg X, \reg X)$ for $n = 1,2$ in other ways.

Noncommutative deformation theory has applications to representation
theory. In Laudal \cite{laud02}, it was shown that noncommutative
deformations of modules are closely related to iterated extensions in
module categories, and we used this result to study finite length
categories of modules in Eriksen \cite{erik04}. These methods work in
any Abelian $k$-category with a reasonable noncommutative deformation
theory.

\emph{Acknowledgments}. The author wishes to thank O.A. Laudal, A.
Siqveland and R. Ile for interesting discussions while preparing
this paper, and the Mittag-Leffler Institute for its hospitality
during the \emph{Noncommutative geometry 2003/04} program.

\section{Functors of noncommutative Artin rings} \label{s:functors}

Let $k$ be an algebraically closed field. We shall define the
category $\acat p$ of \emph{$p$-pointed noncommutative Artin rings}
for any integer $p \ge 1$. For expository purposes, we first define
$\Acat p$, the category of $p$-pointed algebras. An object of $\Acat
p$ is an associative ring $R$, together with structural ring
homomorphisms $f: k^p \to R$ and $g: R \to k^p$ such that $g \circ f
= \id$, and a morphism $u: (R,f,g) \to (R',f',g')$ in $\Acat p$ is
ring homomorphism $u: R \to R'$ such that $u \circ f = f'$ and $g'
\circ u = g$.

We denote by $I = I(R)$ the ideal $I = \ker(g)$ for any $R \in \Acat
p$, and call it the radical ideal of $R$. The category $\acat p$ is the
full subcategory of $\Acat p$ consisting of objects $R \in \Acat p$
such that $R$ is Artinian and (separated) complete in the $I$-adic
topology. For any integer $n \ge 1$, $\acat p(n)$ is the full
subcategory of $\acat p$ consisting of objects $R \in \acat p$ such
that $I^n = 0$. The pro-category $\cacat p$ is the full subcategory of
$\Acat p$ consisting of objects $R \in \Acat p$ such that $R_n = R/I^n$
is Artinian for all $n \ge 1$ and $R$ is (separated) complete in the
$I$-adic topology. It follows that $\acat p \subseteq \cacat p$.

For any $R \in \Acat p$, $R \in \acat p$ if and only if $\dim_k R$ is
finite and $I = I(R)$ is nilpotent. If this is the case, then $I$ is
the Jacobson radical of $R$, and there are $p$ isomorphism classes of
simple left $R$-modules, all of dimension $1$ over $k$.

For any object $R \in \Acat p$, we write $e_1, \dots, e_p$ for the
indecomposable idempotents in $k^p$ and $R_{ij} = e_i R e_j$. Note
that $R$ is a matrix ring in the sense that there is a $k$-linear
isomorphism $R \cong ( R_{ij} ) = \osum R_{ij}$, and multiplication
in $R$ corresponds to matric multiplication in $( R_{ij} )$. In what
follows, we shall denote the direct sum of any family $\{ V_{ij}: 1
\le i,j \le p \}$ of $k$-linear vector spaces by $( V_{ij} )$.

We define a \emph{functor of (p-pointed) noncommutative Artin rings}
to be a covariant functor $\df: \acat p \to \sets$ such that
$\df(k^p) = \{ * \}$ is reduced to one element. It follows that
there is a distinguished element $*_R \in \df(R)$ given by $\df(k^p
\to R)(*)$ for any $R \in \acat p$. We see that $*_R \in \df(R)$ is
a lifting of $* \in \df(k^p)$ to $R$, and call it the \emph{trivial
lifting}.

There is a natural extensions of $\df: \acat p \to \sets$ to the
pro-category $\cacat p$, which we denote by $\df: \cacat p \to
\sets$. For any $R \in \cacat p$, it is given by
    \[ \df(R) = {\underleftarrow{\lim}} \, \df(R_n). \]
A pro-couple for $\df: \acat p \to \sets$ is a pair $(H,\xi)$ with
$H \in \cacat p$ and $\xi \in \df(H)$, and a morphism $u: (H,\xi)
\to (H', \xi')$ of pro-couples is a morphism $u: H \to H'$ in
$\cacat p$ such that $\df(u)(\xi) = \xi'$. By Yoneda's lemma, $\xi
\in \df(H)$ corresponds to a morphism $\phi: \mor(H,-) \to \df$ of
functors on $\acat p$. We say that $(H, \xi)$ \emph{pro-represents}
$\df$ if $\phi$ is an isomorphism of functors on $\acat p$, and that
$(H,\xi)$ is a \emph{pro-representing hull} of $\df$ if $\phi:
\mor(H,-) \to \df$ is a smooth morphism of functors on $\acat p$
that induces an isomorphism of functors on $\acat p(2)$ by
restriction.

\begin{lem}
Let $\df: \acat p \to \sets$ be a functor of noncommutative Artin
rings. If $\df$ has a pro-representing hull, then it is unique up to
a (non-canonical) isomorphism of pro-couples.
\end{lem}
\begin{proof}
Let $(H,\xi)$ and $(H',\xi')$ be pro-representing hulls of $\df$, and
let $\phi, \phi'$ be the corresponding morphisms of functors on $\acat
p$. By the smoothness of $\phi, \phi'$, it follows that $\phi_{H'}$ and
$\phi'_H$ are surjective. Hence there are morphisms $u: (H,\xi) \to
(H', \xi')$ and $v: (H', \xi') \to (H,\xi)$ of pro-couples. Restriction
to $\acat p(2)$ gives morphisms $u_2: (H_2, \xi_2) \to (H'_2,\xi'_2)$
and $v_2: (H'_2, \xi'_2) \to (H_2, \xi_2)$. But both $(H_2, \xi_2)$ and
$(H'_2, \xi'_2)$ represent the restriction of $\df$ to $\acat p(2)$, so
$u_2$ and $v_2$ are mutual inverses. Let us write $\gr_n(R) = I(R)^n /
I(R)^{n+1}$ for all $R \in \cacat p$ and all $n \ge 1$. By the above
argument, it follows that $\gr_1(u)$ and $\gr_1(v)$ are mutual
inverses. In particular, $\gr_1(u \circ v) = \gr_1(u) \circ \gr_1(v)$
is surjective. This implies that $\gr_n(u \circ v)$ is a surjective
endomorphism of the finite dimensional vector space $\gr_n(H')$ for all
$n \ge 1$, and hence an automorphism of $\gr_n(H')$ for all $n \ge 1$.
So $u \circ v$ is an automorphism, and by a symmetric argument, $v
\circ u$ is an automorphism as well. It follows that $u$ and $v$ are
isomorphisms of pro-couples.
\end{proof}

For $1 \le i,j \le p$, let $k^p[\epsilon_{ij}]$ be the object in
$\acat p(2)$ defined by $k^p[\epsilon_{ij}] = k^p + k \cdot
\epsilon_{ij}$, with $\epsilon_{ij} = e_i \, \epsilon_{ij} \, e_j$
and $\epsilon_{ij}^2 = 0$. We define the \emph{tangent space} of
$\df: \acat p \to \sets$ to be $t(\df) = ( t(\df)_{ij} )$ with
$t(\df)_{ij} = \df(k^p[\epsilon_{ij}])$ for $1 \le i,j \le p$. Note
that if $(H,\xi)$ is a pro-representing hull for $\df$, then $\xi$
induces a bijection $t(\df)_{ij} \cong I(H)_{ij} / I(H)^2_{ij}$. In
particular, $t(\df)_{ij}$ has a canonical $k$-linear structure in
this case.

A \emph{small surjection} in $\acat p$ is a surjective morphism $u: R
\to S$ in $\acat p$ such that $KI = IK = 0$, where $I = I(R)$ and $K =
\ker(u)$. Given a functor $\df: \acat p \to \sets$ of noncommutative
Artin rings, a small lifting situation for $\df$ is defined by a small
surjection $u: R \to S$ in $\acat p$ and an element $\xi_S \in \df(S)$.
In order to study the existence of, and ultimately construct, a
pro-representing hull $H$ for $\df$, we are led to consider the
possible liftings of $\xi_S$ to $R$ in small lifting situations.

Let $\{ \ch^n_{ij}: 1 \le i,j \le p \}$ be a family of vector spaces
over $k$ for $n = 1,2$. We say that a functor $\df: \acat p \to
\sets$ of noncommutative Artin rings has an \emph{obstruction
theory} with cohomology $(\ch^n_{ij})$ if the following conditions
hold:
\begin{enumerate}
\item For any small lifting situation, given by a small surjection
$u: R \to S$ in $\acat p$ with kernel $K = \ker(u)$ and an element
$\xi_S \in \df(S)$, we have:
\begin{enumerate}
\item There exists a canonical obstruction $o(u,\xi_S) \in
( \ch^2_{ij} \otimes_k K_{ij} )$ such that $o(u,\xi_S) = 0$ if and
only if there exists a lifting of $\xi_S$ to $R$,
\item If $o(u,\xi_S)=0$, there is an transitive and effective
action of $( \ch^1_{ij} \otimes_k K_{ij} )$ on the set of liftings
of $\xi_S$ to $R$.
\end{enumerate}
\item Let $u_i: R_i \to S_i$ be a small surjection with kernel $K_i
= \ker(u_i)$ and let $\xi_i \in \df(S_i)$ for $i = 1,2$. If $\alpha:
R_1 \to R_2$ and $\beta: S_1 \to S_2$ are morphisms in $\acat p$
such that $u_2 \circ \alpha = \beta \circ u_1$ and
$\df(\beta)(\xi_1) = \xi_2$, then $\alpha^*(o(u_1,\xi_1)) = o(u_2,
\xi_2)$, where $\alpha^*: ( \ch^2_{ij} \otimes_k K_{1,ij} ) \to (
\ch^2_{ij} \otimes_k K_{2,ij} )$ is the natural map induced by
$\alpha$.
\end{enumerate}
Moreover, if $\ch^n_{ij}$ has finite $k$-dimension for $1 \le i,j
\le p, \; n = 1,2$, then we say that $\df$ has an obstruction theory
with \emph{finite dimensional cohomology} $(\ch^n_{ij})$.

In the rest of this section, we shall assume that $\df: \acat p \to
\sets$ is a functor of noncommutative Artin rings that has an
obstruction theory with finite dimensional cohomology $( \ch^n_{ij}
)$. Note that for any object $R \in \acat p(2)$, the morphism $R \to
k^p$ is a small surjection. This implies that there is a canonical
set-theoretical bijection
\[ ( \ch^1_{ij} \otimes_k I(R)_{ij} ) \cong \df(R), \]
given by the trivial lifting $*_R \in \df(R)$. In particular, there
is a set-theoretical bijection between $H^1_{ij}$ and $t(\df)_{ij} =
\df(k^p[\epsilon_{ij}])$ for $1 \le i,j \le p$.

We define $\mathit{T}^n$ to be the free, formal matrix ring in
$\cacat p$ generated by the $k$-linear vector spaces $ \{
(\ch^n_{ij})^*: 1 \le i,j \le p \}$ for $n=1,2$, where
$(\ch^n_{ij})^* = \hmm_k(\ch^n_{ij}, k)$. For any $R \in \acat
p(2)$, we have natural isomorphisms
\[ \df(R) \cong ( \ch^1_{ij} \otimes_k I(R)_{ij} ) \cong (
\hmm_k((\ch^1_{ij})^*, I(R)_{ij}) ) \cong \mor(\infdef_2, R), \]
where $\infdef_2 = \infdef / I(\infdef)^2$. It follows that there is
an isomorphism $\mor(\infdef_2,-) \to \df$ of functors on $\acat
p(2)$, i.e. the restriction of $\df$ to $\acat p(2)$ is represented
by $(\infdef_2, \xi_2)$ for some $\xi_2 \in \df(\infdef_2)$.

\begin{thm} \label{t:obstrmor}
Let $\df: \acat p \to \sets$ be a functor of noncommutative Artin
rings. If $\df$ has an obstruction theory with finite dimensional
cohomology, then there is an obstruction morphism $o: \obstrdef \to
\infdef$ in $\cacat p$ such that $H = \infdef \widehat
\otimes_{\obstrdef} k^p$ is a pro-representing hull of $\df$.
\end{thm}
\begin{proof}
Let us write $I = I(\infdef)$, $\infdef_n = \infdef / I^n$ and $t_n:
\infdef_{n+1} \to \infdef_n$ for the natural morphism for all $n \ge
1$. Let $a_2 = I^2$ and $H_2 = \infdef / a_2 = \infdef_2$, then the
restriction of $\df$ to $\acat p(2)$ is represented by $(H_2,\xi_2)$
and $H_2 \cong \infdef_2 \otimes_{\obstrdef} k^p$. Using $o_2$ and
$\xi_2$ as a starting point, we shall construct $o_{n+1}$ and
$\xi_{n+1}$ for $n \ge 2$ inductively. So let $n \ge 2$, and assume
that the morphism $o_n: \obstrdef \to \infdef_n$ and the deformation
$\xi_n \in \df(H_n)$ are given, with $H_n = \infdef_n
\otimes_{\obstrdef} k^p$. We may assume that $t_{n-1} \circ o_n =
o_{n-1}$ and that $\xi_n$ is a lifting of $\xi_{n-1}$.

Let us first construct the morphism $o_{n+1}: \obstrdef \to
\infdef_{n+1}$. We define $a'_n$ to be the ideal in $\infdef_n$
generated by $o_n(I(\obstrdef))$. Then $a'_n = a_n/I^n$ for an
ideal $a_n \subseteq \infdef$ with $I^n \subseteq a_n$, and $H_n
\cong \infdef / a_n$. Let $b_n = I a_n + a_n I$, then we obtain
the following commutative diagram:
\[
\xymatrix{ {\obstrdef} \ar[rd]^{o_n} & {\infdef_{n+1}} \ar[d]
\ar[r] & {\infdef/b_n}
\ar[d] \\
& {\infdef_n} \ar[r] & H_n=\infdef/a_n, }
\]
There is an obstruction $o'_{n+1}=o(\infdef/b_n \to H_n,\xi_n)$ for
lifting $\xi_n$ to $\infdef/b_n$ since $\infdef/b_n \to \infdef/a_n$
is a small surjection, hence a morphism $o'_{n+1}: \obstrdef \to
\infdef/b_n$. Let $a''_{n+1}$ be the ideal in $\infdef/b_n$
generated by $o'_{n+1}(I(\obstrdef))$. Then $a''_{n+1}=a_{n+1}/b_n$
for an ideal $a_{n+1} \subseteq \infdef$ with $b_n \subseteq a_{n+1}
\subseteq a_n$. Let $H_{n+1}=\infdef/a_{n+1}$, then we obtain the
following commutative diagram:
\[
\xymatrix{ {\obstrdef} \ar[rd]_{o_n} \ar@/^1pc/[rr]^{o'_{n+1}} &
{\infdef_{n+1}} \ar[d] \ar[r] & {\infdef/b_n} \ar[r] \ar[d] &
H_{n+1}=
\infdef/a_{n+1} \ar[dl] \\
& {\infdef_n} \ar[r] & H_n=\infdef/a_n & }
\]
By the choice of $a_{n+1}$, the obstruction for lifting $\xi_n$ to
$H_{n+1}$ is zero. We can therefore find a lifting $\xi_{n+1} \in
\df(H_{n+1})$ of $\xi_n$ to $H_{n+1}$.

We claim that there is a morphism $o_{n+1}:\obstrdef \to
\infdef_{n+1}$ that commutes with $o'_{n+1}$ and $o_n$. Note that
$a_{n-1} = I^{n-1} + a_n$ since $t_{n-1} \circ o_n = o_{n-1}$. For
simplicity, we write $O(K)=(\hmm_k(\gr_1(\obstrdef)_{ij},K_{ij}))$
for any family $K = (K_{ij})$ of vector spaces over $k$. The
following diagram of $k$-vector spaces is commutative with exact
columns:
\[
\xymatrix{
{0} \ar[d] & {0} \ar[d] \\
{O(b_n/I^{n+1})} \ar[d] \ar[r]^{j_n} & {O(b_{n-1}/I^n)} \ar[d] \\
{O(a_n/I^{n+1})} \ar[d]^{r_{n+1}} \ar[r]^{k_n} & {O(a_{n-1}/I^n)}
\ar[d]^{r_n} \\
{O(a_n/b_n)} \ar[d] \ar[r]^{l_n} & {O(a_{n-1}/b_{n-1})} \ar[d] \\
{0} & {0} }
\]
We may consider $o_n$ as an element in $O(a_{n-1}/I^n)$, and
$o'_{n+1} \in O(a_n/b_n)$. Since $o'_n$ commutes with $o'_{n+1}$ and
$o_n$, we get $l_n(o'_{n+1}) = r_n(o_n)$. To prove the claim, it is
enough to find an element $o_{n+1} \in O(a_n/I^{n+1})$ such that
$r_{n+1}(o_{n+1}) = o'_{n+1}$ and $k_n(o_{n+1}) = o_n$. Since
$o_n(I(\obstrdef)) \subseteq a_n$, there is an element
$\overline{o}_{n+1} \in O(a_n/I^{n+1})$ such that
$k_n(\overline{o}_{n+1})=o_n$. But $a_{n-1}=a_n + I^{n-1}$ implies
that $j_n$ is surjective, so the claim follows from the snake lemma.
In particular, $\infdef_{n+1} \otimes_{\obstrdef} k^p \cong H_{n+1}$
when the tensor product is taken over $o_{n+1}$.

By induction, we find a morphism $o_n: \obstrdef \to \infdef_n$ and
an element $\xi_n \in \df(H_n)$ for all integers $n \ge 1$, with
$H_n = \infdef_n \otimes_{\obstrdef} k^p$. Using the universal
property of the projective limit, we obtain a morphism $o: \obstrdef
\to \infdef$ in $\cacat p$ and an element $\xi \in \df(H)$, with $H
= \infdef \hat{\otimes}_{\obstrdef} k^p$. We claim that $(H, \xi)$
is a pro-representable hull for $\df$.

Clearly, it is enough to prove that $(H_n, \xi_n)$ is a
pro-representing hull for the restriction of $\df$ to $\acat p(n)$
for all $n\ge 3$. So let $\phi_n: \mor(H_n,-) \to \df$ be the
morphism of functors on $\acat p(n)$ corresponding to $\xi_n$ for
some $n \ge 3$. We shall prove that $\phi_n$ is a smooth morphism.
Let $u: R \to S$ be a small surjection in $\acat p(n)$ with kernel
$K$, let $E_R \in \df(R)$ and $v \in \mor(H_n,S)$ be elements such
that $\df(u)(E_R) = \df(v)(\xi_n) = E_S$, and consider the following
commutative diagram:
\[
\xymatrix{
{\infdef} \ar[r] \ar[rd] & {H_{n+1}} \ar[d] &
{R} \ar[d]^u \\
& {H_n} \ar[r]_v & {S} }
\]
We can find a morphism $v': \infdef \to R$ that makes the diagram
commutative. This implies that $v'(a_n) \subseteq K$, and since $u$ is
small, that $v'(b_n) = 0$. But the induced map $\infdef/b_n \to R$ maps
the obstruction $o'_{n+1}$ to $o(u, E_S) = 0$. It follows that
$v'(a_{n+1}) = 0$, hence $v'$ induces a morphism $v': H_{n+1} \to R$
making the diagram commutative. Since $v'(I(H_{n+1})^n)=0$, we may
consider $v'$ as a map from $H_{n+1}/I(H_{n+1})^n \cong H_n$. This
proves that there is a morphism $v':H_n \to R$ such that $u \circ v' =
v$.

Let $E'_R = \df(v')(\xi_n)$, then $E'_R$ is a lifting of $E_S$ to
$R$, and the difference between $E_R$ and $E'_R$ is given by an
element $d \in (\ch^1_{ij} \otimes_k K_{ij}) =
(\hmm_k(\gr_1(\infdef)_{ij},K_{ij}))$. Let $v'': \infdef \to R$ be
the morphism given by $v''(x_{ij}(l)) = v'(x_{ij}(l)) +
d(\overline{x_{ij}(l)})$, where $\{ x_{ij}(l): 1 \le l \le d_{ij}
\}$ is a basis for $\ch^1_{ij}$ for $1 \le i,j \le p$. Since
$a_{n+1} \subseteq I^2$ and $u$ is small, $v''(a_{n+1}) \subseteq
v'(a_{n+1}) + I(R)K + KI(R) + K^2 = v'(a_{n+1}) = 0$. This implies
that $v''$ induces a morphism $v'': H_n \to R$. By construction, $u
\circ v'' = u \circ v' = v$ and $\df(v'') (\xi_n) = E_R$, and this
proves that $\phi_n$ is smooth.
\end{proof}

We remark that a more general version of theorem \ref{t:obstrmor}
can be proved if $\df$ has an obstruction theory with cohomology $(
\ch^n_{ij} )$ and $\ch^1_{ij}$ has a countable $k$-basis for $1 \le
i,j \le p$, using the methods of Laudal \cite{laud79}.

\begin{cor} \label{c:tan}
Let $\df: \acat p \to \sets$ be a functor of noncommutative Artin
rings. If $\df$ has an obstruction theory with finite dimensional
cohomology $( \ch^n_{ij} )$, then there is a $k$-linear isomorphism
$\ch^1_{ij} \cong t(\df)_{ij}$ for $1 \le i,j \le p$.
\end{cor}

\begin{cor} \label{c:massey}
Let $\df: \acat p \to \sets$ be a functor of noncommutative Artin
rings. If $\df$ has an obstruction theory with finite dimensional
cohomology $( \ch^n_{ij} )$, then the pro-representing hull
$\hull(\df)$ is completely determined by the $k$-linear spaces
$\ch^n_{ij}$ for $1 \le i,j \le p, \; n = 1,2$, together with some
generalized Massey products on them.
\end{cor}

These results are natural generalizations of similar results for
functors of commutative Artin rings. As in the commutative case, the
generalized Massey product structure on $\ch^n_{ij}$ can be considered
as the $k$-linear dual of the obstruction morphism $o: \obstrdef \to
\infdef$, see Laudal \cite{laud79}, \cite{laud86}. If $\df$ is
obstructed, i.e. $o(I(\obstrdef)) \neq 0$, then it is a non-trivial
task to compute $\hull(\df)$ using generalized Massey products.

Let $\df: \acat p \to \sets$ be a functor of noncommutative Artin
rings. For $1 \le i \le p$, we write $\df_i: \acat 1 \to \sets$ for the
restriction of $\df$ to $\acat 1$ using the $i$'th natural inclusion of
categories $\acat 1 \hookrightarrow \acat p$, and $\df^c_i: \lcat \to
\sets$ for the restriction of $\df_i$ to $\lcat$, where $\lcat$ is the
full subcategory of $\acat 1$ consisting of commutative algebras. For
any associative ring $R$, we define the commutativization of $R$ to be
the quotient ring $R^c = R/I^c(R)$, where $I^c(R)$ is the ideal in $R$
generated by the set of commutators $\{ ab-ba: a,b \in R \}$.

\begin{prop} \label{p:commhull}
Let $\df: \acat p \to \sets$ be a functor of noncommutative Artin
rings. If $\df$ has an obstruction theory with finite dimensional
cohomology $(\ch^n_{ij})$, then $\df^c_i$ has a commutative hull
$\hull(\df^c_i) \cong \hull(\df_i)^c$ for $1 \le i \le p$, and
$\hull(\df)^c \cong \oplus \hull(\df^c_i)$.
\end{prop}
\begin{proof}
Clearly, the functor $\df_i$ has an obstruction theory with finite
dimensional cohomology $\ch^n_{ii}$ for $1 \le i \le p$, and therefore
a pro-representing hull $\hull(\df_i)$ that is determined by an
obstruction morphism $o_i: \obstrdef \to \infdef$. Similarly, the
functor $\df^c_i$ has an obstruction theory with finite dimensional
cohomology $\ch^n_{ii}$ for $1 \le i \le p$, and therefore a
pro-representing hull $\hull(\df^c_i)$ that is determined by an
obstruction morphism $o^c_i: (\obstrdef)^c \to (\infdef)^c$. These
morphisms are defined by obstructions in small lifting situations, so
it follows from the functorial nature of the obstructions that $o_i$
and $o^c_i$ are compatible. Hence $\hull(\df^c_i) \cong \hull(\df_i)^c$
for $1 \le i \le p$. For the second part, note that $\hull(\df)^c_{ij}
= 0$ whenever $i \neq j$. In fact, for any $x_{ij} \in \hull(\df)_{ij}$
with $i \neq j$, the commutator $[e_i, x_{ij}] = e_i x_{ij} - x_{ij}
e_i = x_{ij}$ is zero in $\hull(\df)^c$. This implies that
$\hull(\df)^c = \oplus \hull(\df)^c_{ii} \cong \oplus \hull(\df_i)^c
\cong \oplus \hull(\df^c_i)$.
\end{proof}

\section{Noncommutative deformation functors} \label{s:defm}

Let $k$ be an algebraically closed field, and let $\cat A$ be any
Abelian $k$-category. For any object $X \in \cat A$, we recall the
definition of $\defm{X}^{\cat A}: \lcat \to \sets$, the commutative
deformation functor of $X$ in $\cat A$, which is a functor of Artin
rings, and discuss how to generalize this definition to
noncommutative deformation functors $\defm{X}^{\cat A}: \acat p \to
\sets$ of a family $X = \{ X_1, \dots, X_p \}$ in the category $\cat
A$.

Let $R$ be any object in $\acat p$. We consider the category $\cat
A_R$ of $R$-objects in $\cat A$, i.e. the category with objects
$(X,\phi)$, where $X$ is an object of $\cat A$ and $\phi: R \to
\mor_{\cat A}(X,X)$ is a $k$-algebra homomorphism, and with
morphisms $f: (X,\phi) \to (X',\phi')$, where $f:X \to X'$ is a
morphism in $\cat A$ such that $f \circ \phi(r) = \phi'(r) \circ f$
for all $r \in R$. Clearly, $\cat A_R$ is an Abelian $R$-category.

Let $R$ be any object of $\acat p$, let $\cat B$ be any Abelian
$R$-category, and let $\smcat R$ be the category of finitely
generated left $R$-modules. For each object $Y \in \cat B$, there is
a unique finite colimit preserving functor $Y \otimes_R -: \smcat R
\to \cat B$ that maps $R$ to $Y$, given in the following way: If $M
= \coker(f)$, where $f: R^m \to R^n$ is a homomorphism of left
$R$-modules, then $Y \otimes_R M = \coker(F)$, where $F: Y^m \to
Y^n$ is the morphism in $\cat B$ induced by $f$ and the $R$-linear
structure on $\cat B$. We say that an object $Y \in \cat B$ is
\emph{$R$-flat} if $Y \otimes_R -$ is exact. It is clear that any
morphism $u: R \to S$ in $\acat p$ induces a functor $- \otimes_R S:
\cat A_R \to \cat A_S$.

Given an object $X \in \cat A$, the flat deformation functor
$\defm{X`}^{\cat A}: \lcat \to \sets$ is given in the following way:
For any object $R \in \lcat$, we define a lifting of $X$ to $R$ to
be an object $X_R \in \cat A_R$ which is $R$-flat, together with an
isomorphism $\eta: X_R \otimes_R k \to X$ in $\cat A$, and we say
that two liftings $(X_R, \eta)$ and $(X'_R, \eta')$ are equivalent
if there is an isomorphism $\tau: X_R \to X'_R$ in $\cat A_R$ such
that $\eta' \circ (\tau \otimes_R \id) = \eta$. Let $\defm{X}^{\cat
A}(R)$ be the set of equivalence classes of liftings of $X$ to $R$.
Then $\defm{X}^{\cat A}: \lcat \to \sets$ is a functor of Artin
rings, called the commutative deformation functor of $X$ in $\cat
A$.

When $\cat A = \mcat A$, the category of left modules over an
associative $k$-algebra $A$, we remark that the category $\cat A_R$
is the category of $A$-$R$ bimodules on which $k$ acts centrally,
and the tensor product defined above is the usual tensor product
over $R$. It follows that the usual deformation functor $\defm M$ of
a left $A$-module $M$ coincides with the deformation functor defined
above.

Given a finite family $X = \{ X_1, \dots, X_p \}$ of objects in
$\cat A$, we would like to define a noncommutative deformation
functor $\defm{X}^{\cat A}: \acat p \to \sets$ of the family $X$ in
$\cat A$. When $\cat A = \mcat A$, the category of left modules over
an associative $k$-algebra $A$, such a deformation functor was
defined in Laudal \cite{laud02}. The idea is to replace the
condition that $M_R$ is a flat right $R$-module with the
\emph{matric freeness} condition that
\begin{equation}
\label{e:mfree} M_R \cong ( M_i \otimes_k R_{ij} )
\end{equation}
as right $R$-modules. This is reasonable, since an $R$-module is
flat if and only if (\ref{e:mfree}) holds when $R \in \lcat$ or $R
\in \acat 1$, see Bourbaki \cite{bour-ac61i}, Corollary II.3.2.
However, it is not clear whether an $R$-module is flat if and only
if (\ref{e:mfree}) holds when $R \in \acat p$ for $p \ge 2$.

We choose to define noncommutative deformation functors using
Laudal's matric freeness condition rather than flatness. However, it
is not completely clear how to do this for an arbitrary Abelian
$k$-category. We shall therefore restrict our attention to
categories of sheaves and presheaves of modules.

Let $\ccat$ be a small category, let $\sa$ be a presheaf of
$k$-algebras on $\ccat$, and assume that $\cat A$ is an Abelian
$k$-category and a full subcategory of the category $\presh(\ccat,
\sa)$ of presheaves of left $\sa$-modules on $\ccat$. Then there is a
forgetful functor $\pi_c: \cat A \to \mcat k$ for each object $c \in
\ccat$ and an induced forgetful functor $\pi_c^R: \cat A_R \to \mcat R$
for each object $c \in \ccat$ and each $R \in \acat p$. We say that an
object $X_R \in \cat A_R$ is \emph{$R$-free} if $\pi^R_c(X_R) \cong (
\pi_c(X_i) \otimes_k R_{ij} )$ in $\mcat R$ for all objects $c \in
\ccat$.

Let $\cat A$ be an Abelian $k$-category, and let $X = \{ X_1, \dots,
X_p \}$ be a finite family of objects in $\cat A$. If $\cat A$ is a
full subcategory of $\presh(\ccat,\sa)$, we define the
noncommutative deformation functor $\defm{X}^{\cat A}: \acat p \to
\sets$ in the following way: A lifting of $X$ to $R$ is an object
$X_R$ in $\cat A_R$ that is R-free, together with isomorphisms
$\eta_i: X_R \otimes_R k_i \to X_i$ in $\cat A$ for $1 \le i \le p$,
and two liftings $(X_R, \eta_i)$ and $(X'_R, \eta'_i)$ are
equivalent if there is an isomorphism $\tau: X_R \to X'_R$ in $\cat
A_R$ such that $\eta'_i \circ (\tau \otimes_R \id) = \eta_i$ for $1
\le i \le p$. Let $\defm{X}^{\cat A}(R)$ be the set of equivalence
classes of liftings of $X$ to $R$. Then $\defm{X}^{\cat A}$ is a
functor of noncommutative Artin rings, the noncommutative
deformation functor of the family $X$ in $\cat A$. When the category
$\cat A$ is understood from the context, we often write $\defm X$
for $\defm{X}^{\cat A}$.

\section{Deformations of presheaves of modules} \label{s:defm-presh}

Let $k$ be an algebraically closed field. The category
$\presh(\ccat, \sa)$ of presheaves of left $\sa$-modules on $\ccat$
is an Abelian $k$-category for any small category $\ccat$ and any
presheaf $\sa$ of associative $k$-algebras on $\ccat$, and we shall
consider deformations in this category. To fix notations, a presheaf
on $\ccat$ is always covariant in this paper.

For any finite family $\sfam = \{ \sh F_1, \dots, \sh F_p \}$ of
presheaves of left $\sa$-modules on $\ccat$, we consider the
noncommutative deformation functor $\defs: \acat p \to \sets$,
defined by $\defs = \defs^{\cat A}$ with $\cat A = \presh(\ccat,
\sa)$. We shall describe this functor in concrete terms.

Let $R \in \acat p$, and consider a lifting $\sh F_R$ of the family
$\sfam$ to $R$. Without loss of generality, we may assume that $\sh
F_R(c) = ( \sh F_i(c) \otimes_k R_{ij})$ with the natural right
$R$-module structure for all $c \in \ccat$. To describe the lifting
completely, we must specify the left action of $\sa(c)$ on $(\sh
F_i(c) \otimes_k R_{ij})$ for any object $c \in \ccat$, and the
restriction map $\sh F_R(\phi): ( \sh F_i(c) \otimes_k R_{ij}) \to (
\sh F_i(c') \otimes_k R_{ij})$ for any morphism $\phi: c \to c'$ in
$\ccat$. It is enough to specify the action of $a \in \sa(c)$ on
elements of the form $f_j \otimes e_j$ in $\sh F_j(c) \otimes_k
R_{jj}$, and we must have
\begin{equation} \label{e:def-la}
a (f_j \otimes e_j) = (af_j) \otimes e_j + \sum f'_i \otimes r_{ij}
\end{equation}
with $f'_i \in \sh F_i(c), \; r_{ij} \in I(R)_{ij}$ for all objects
$c \in \ccat$. Similarly, it is enough to specify the restriction
map $\sh F_R(\phi)$ on elements of the form $f_j \otimes e_j$ in
$\sh F_j(c) \otimes_k R_{jj}$, and we must have
\begin{equation} \label{e:def-res}
\sh F_R(\phi)(f_j \otimes e_j) = \sh F_j(\phi)(f_j) \otimes e_j +
\sum f'_i \otimes r_{ij}
\end{equation}
with $f'_i \in \sh F_i(c'), \; r_{ij} \in I(R)_{ij}$ for all
morphisms $\phi: c \to c'$ in $\ccat$.

Let $Q^R(c,c') = ( \hmm_k(\sh F_j(c), \sh F_i(c') \otimes_k R_{ij})
)$ for all objects $c,c' \in \ccat$, and write $Q^R(c) = Q^R(c,c)$.
There is a natural product $Q^R(c',c'') \otimes_k Q^R(c,c') \to
Q^R(c,c'')$ for all objects $c,c',c'' \in \ccat$, given by
composition of maps and multiplication in $R$, such that $Q^R(c)$ is
an associative $k$-algebra and $Q^R(c,c')$ is an $Q^R(c')$-$Q^R(c)$
bimodule in a natural way.

\begin{lem} \label{l:defm-char}
For any $R \in \acat p$, there is a bijective correspondence between
the following data, up to equivalence, and $\defs(R)$:
\begin{enumerate}
    \item For any $c \in \ccat$, a $k$-algebra homomorphism $L(c):
    \sa(c) \to Q^R(c)$ that satisfies equation (\ref{e:def-la}),
    \item For any morphism $\phi: c \to c'$ in $\ccat$, an element
    $L(\phi) \in Q^R(c,c')$ that satisfies equation (\ref{e:def-res})
    and $L(\phi) L(c) = L(c') L(\phi)$,
    \item We have $L(\id) = \id$ and $L(\phi') L(\phi) = L(\phi'
    \circ \phi)$ for all morphisms $\phi: c \to c'$ and $\phi': c'
    \to c''$ in $\ccat$.
\end{enumerate}
\end{lem}

For any $R \in \acat p$, we may for any $c \in \ccat$ consider $L(c):
\sa(c) \to Q^R(c)$ given by $L(c)(a)(f_j) = af_j \otimes e_j$ for all
$a \in \sa(c), \; f_j \in \sh F_j(c)$, and for any morphism $\phi: c
\to c'$ in $\ccat$ consider $L(\phi) \in Q^R(c,c')$ given by
$L(\phi)(f_j) = \sh F_j(\phi)(f_j) \otimes e_j$ for all $f_j \in \sh
F_j(c)$. These data correspond to the trivial deformation $*_R \in
\defs(R)$.

\begin{lem} \label{l:tg-ext}
There is a bijection $\phi_{ij}: t(\defs)_{ij} \to \ext^1_{\sa}(\sh
F_j, \sh F_i)$ for $1 \le i,j \le p$ that maps trivial deformations to
split extensions.
\end{lem}
\begin{proof}
Let $R = k^p[\epsilon_{ij}]$, and let $\sh F_R$ be a lifting of $\sfam$
to $R$. We consider the $j$'th column $\sh F^j_R$ of $\sh F_R$, given
by $c \mapsto \sh F_R(c) e_j$, which is a sub-presheaf of $\sh F_R$ of
left $\sa$-modules on $\ccat$ since $\sh F^j_R(c)$ is invariant under
$L(c)$ and $L(\phi)$ for any $c \in \ccat$ and any $\phi: c \to c'$ in
$\ccat$. Moreover, there is a natural exact sequence $0 \to \sh F_i \to
\sh F_R^j \to \sh F_j \to 0$ in $\presh(\ccat, \sa)$, since $\sh
F^j_R(c) = \sh F_j(c) \oplus \sh F_i(c) \epsilon_{ij}$ for all $c \in
\ccat$. Clearly, equivalent liftings of $\sfam$ to $R$ give equivalent
extensions in $\presh(\ccat, \sa)$, so $\sh F_R \mapsto \sh F^j_R$
defines a map $\phi_{ij}: t(\defs)_{ij} \to \ext^1_{\sa}(\sh F_j, \sh
F_i)$ that maps trivial deformations to split extensions. To construct
an inverse of $\phi_{ij}$, we consider an extension $\sh E$ of $\sh
F_j$ by $\sh F_i$ in $\presh(\ccat, \sa)$, and let $\sh F_R = \sh E
\oplus \sh F_1 \dots \oplus \widehat{\sh F_i} \oplus \dots
\widehat{\oplus \sh F_j} \oplus \dots \oplus \sh F_p$. Then $\sh F_R$
is a presheaf of left $\sa$-modules on $\ccat$, and it is easy to see
that it defines a lifting of $\sfam$ to $R$ since $\sh E(c) \cong \sh
F_i(c) \oplus \sh F_j(c)$ as $k$-linear vector spaces for any $c \in
\ccat$. It follows that the assignment $\sh E \mapsto \sh F_R$ defines
an inverse of $\phi_{ij}$.
\end{proof}

\section{Global Hochschild cohomology} \label{s:coh}

Let $k$ be an algebraically closed field, let $\ccat$ be a small
category, and let $\sa$ be a presheaf of associative $k$-algebras on
$\ccat$. For any presheaves $\sh F, \sh G$ of left $\sa$-modules on
$\ccat$, we may consider the Hochschild complex $\hc^{\cpd}(\sa(c),
\hmm_k(\sh F(c), \sh G(c)))$ of $\sa(c)$ with values in the bimodule
$\hmm_k(\sh F(c),\sh G(c))$ for any object $c \in \ccat$.
Unfortunately, this construction is not functorial in $c$. In this
section, we consider a variation of this construction that is
functorial, and use this to define a global Hochschild cohomology
theory.

Let $\morc$ denote the category of morphisms in $\ccat$ defined in the
following way: An object in $\morc$ is a morphism in $\ccat$, and given
objects $f: c \to c'$ and $g: d \to d'$ in $\morc$, a morphism
$(\alpha, \beta): f \to g$ in $\morc$ is a couple of morphisms $\alpha:
d \to c$ and $\beta: c' \to d'$ in $\ccat$ such that $\beta f \alpha =
g$. Clearly, $\morc$ is a small category.

Let $\sh F, \sh G$ be presheaves of left $\sa$-modules on $\ccat$,
and let $\shhom_k(\sh F, \sh G)$ be the presheaf in
$\presh(\morc,\underline k)$ given by $\shhom_k(\sh F, \sh G)(\phi)
= \hmm_k(\sh F(c), \sh G(c'))$ for any morphism $\phi: c \to c'$ in
$\ccat$. We define the \emph{Hochschild complex} of $\sa$ with
values in $\shhom_k(\sh F, \sh G)$ to be the functor
$\hc^{\cpd}(\sa, \shhom_k(\sh F, \sh G)): \morc \to \compl(k)$,
given by
    \[ \hc^p(\sa, \shhom_k(\sh F, \sh G))(\phi) =
    \hmm_k(\otimes^p_k \sa(c), \hmm_k(\sh F(c), \sh G(c'))) \]
for any morphism $\phi: c \to c'$ in $\ccat$ and any integer $p \ge 0$,
with differential given by
\begin{equation*}
\begin{split}
d^p(\phi)(f)(a_1 \otimes \dots \otimes a_{p+1}) & = a_1 \; f(a_2
\otimes \dots \otimes a_{p+1}) \\
& + \sum_{i=1}^p (-1)^i \; f(a_1 \otimes \dots \otimes a_i a_{i+1}
\otimes \dots \otimes a_{p+1}) \\
& + (-1)^{p+1} \; f(a_1 \otimes \dots \otimes a_p) \; a_{p+1}
\end{split}
\end{equation*}
for any $f \in \hc^p(\sa, \shhom_k(\sh F, \sh G))(\phi)$ and any
$a_1, \dots, a_{p+1} \in \sa(c)$. We see that $\hc^{\cpd}(\sa,
\shhom_k(\sh F, \sh G))(\phi)$ is the Hochschild complex of $\sa(c)$
with values in the bimodule $\hmm_k(\sh F(c), \sh G(c'))$ for any
morphism $\phi: c \to c'$ in $\ccat$, and the definition of $\morc$
ensures that $\phi \mapsto \hc^{\cpd}(\sa, \shhom_k(\sh F, \sh
G))(\phi)$ is functorial.

Let us write $\shext^{\cpd}_{\sa}(\sh F, \sh G): \morc \to \mcat k$
for the composition of functors given by $\phi \mapsto
\hc^{\cpd}(\sa, \shhom_k(\sh F, \sh G))(\phi) \mapsto
\ch^{\cpd}(\hc^{\cpd}(\sa, \shhom_k(\sh F, \sh G))(\phi))$. We see
that $\shext^{\cpd}_{\sa}(\sh F, \sh G)(\phi) \cong
\ext^{\cpd}_{\sa(c)}(\sh F(c), \sh G(c'))$ for any morphism $\phi: c
\to c'$ in $\ccat$.

For any functor $G: \morc \to \mcat k$, we may consider the
resolving complex $\dc^{\cpd}(\ccat,G)$ in $\mcat k$ of the
projective limit functor of $G$, see Laudal \cite{laud79}. We recall
that for any integer $p \ge 0$, $\dc^p(\ccat, G)$ is given by
    \[ \dc^p(\ccat,G) = \prod_{c_0 \to \dots \to c_p} G(\phi_p
    \circ \dots \circ \phi_1), \]
where the product is taken over all $p$-tuples $(\phi_1, \dots,
\phi_p)$ of composable morphisms $\phi_i: c_{i-1} \to c_i$ in
$\ccat$, and the differential $d^p: \dc^p(\ccat,G) \to
\dc^{p+1}(\ccat,G)$ is given by
\begin{equation*}
\begin{split}
(d^p g)(\phi_1, \dots, \phi_{p+1}) & =  \; G(\phi_1,\id)(g(\phi_2,
\dots, \phi_{p+1})) \\
& + \sum_{i=1}^p \; (-1)^i \; g(\phi_1, \dots, \phi_{i+1} \circ
\phi_i, \dots, \phi_{p+1}) \\
& + (-1)^{p+1} \; G(\id,\phi_{p+1})(g(\phi_1, \dots, \phi_p))
\end{split}
\end{equation*}
for all $g \in \dc^p(\ccat,G)$ and for all ($p+1$)-tuples $(\phi_1,
\dots, \phi_{p+1})$ of composable morphisms $\phi_i: c_{i-1} \to c_i$
in $\ccat$. We denote the cohomology of $\dc^{\cpd}(\ccat, G)$ by
$\ch^{\cpd}(\ccat, G)$, and recall the following standard result:

\begin{prop} \label{p:dcmpl}
Let $\ccat$ be a small category. The resolving complex
$\dc^*(\ccat,G)$ has the following properties:
\begin{enumerate}
    \item $\dc^{\cpd}(\ccat,-): \presh(\morc, \underline k) \to
    \compl(k)$ is exact,
    \item $\ch^p(\ccat,G) \cong \underleftarrow{\lim}^{(p)} G$ for
    all $G \in \presh(\morc, \underline k)$ and for all $p \ge 0$.
\end{enumerate}
In particular, $\ch^{\cpd}(\ccat,-): \presh(\morc, \underline k) \to
\mcat k$ is an exact $\delta$-functor.
\end{prop}

For any functor $\cc^{\cpd}: \morc \to \compl(k)$, we may consider the
double complex $\dc^{\cpd \cpd} = \dc^{\cpd}(\ccat, \cc^*)$ of vector
spaces over $k$. Explicitly, have that $\dc^{pq} = \dc^p(\ccat,\cc^q)$
for all integers $p,q$ with $p \ge 0$, that $\dI{pq}: \dc^{pq} \to
\dc^{p+1,q}$ is the differential $d^p$ in $\dc^{\cpd}(\ccat, \cc^q)$,
and that $\dII{pq}: \dc^{pq} \to \dc^{p,q+1}$ is the differential given
by $\dII{pq} = (-1)^p \; \dc^p(\ccat, d^q)$, where $d^q: \cc^q \to
\cc^{q+1}$ is the differential in $\cc^{\cpd}$. Note that if $\cc^q =
0$ for all $q<0$, then $\dc^{\cpd \cpd}$ lies in the first quadrant.

We define the \emph{global Hochschild complex} of $\sa$ with values
in $\shhom_k(\sh F, \sh G)$ on $\ccat$ to be the total complex of
the double complex $\dc^{\cpd \cpd} = \dc^{\cpd}(\ccat,
\hc^{\cpd}(\sa, \shhom_k(\sh F, \sh G)))$, and denote it by
$\hc^{\cpd}(\ccat, \sa, \shhom_k(\sh F, \sh G))$. Moreover, we
define \emph{global Hochschild cohomology} $\hh^{\cpd}(\ccat, \sa,
\shhom_k(\sh F, \sh G))$ of $\sa$ with values in $\shhom_k(\sh F,
\sh G)$ on $\ccat$ to be the cohomology of the global Hochschild
complex $\hc^{\cpd}(\ccat, \sa, \shhom_k(\sh F, \sh G))$.

\begin{prop} \label{p:glh-ss}
There is a spectral sequence for which $E^{pq}_2 = \ch^p(\ccat,
\shext^q_{\sa}(\sh F, \sh G))$, such that $E_{\infty} = \gr
\hh^{\cpd}(\ccat, \sa, \shhom_k(\sh F, \sh G))$, the associated
graded vector space over $k$ with respect to a suitable filtration
of $\hh^{\cpd}(\ccat, \sa, \shhom_k(\sh F, \sh G))$.
\end{prop}

We remark that $\hh^{\cpd}(\ccat, \sa, \shhom_k(\sh F, \sh G))$ can
be calculated in concrete terms in many situations, using the above
spectral sequence. We shall give an example of such a computation in
last section of this paper.

\section{Obstruction theory for presheaves of modules}
\label{s:obstr}

Let $k$ be an algebraically closed field. For any small category
$\ccat$ and any presheaf $\sa$ of associative $k$-algebras on
$\ccat$, we shall construct an obstruction theory for the
noncommutative deformation functor $\defs: \acat p \to \sets$ with
cohomology $( \hh^n(\ccat, \sa, \shhom_k(\sh F_j, \sh F_i)) )$ for
any finite family $\sfam = \{ \sh F_1, \dots, \sh F_p \}$ of
presheaves of left $\sa$-modules on $\ccat$.

\begin{prop} \label{p:obstr}
Let $u: R \to S$ be a small surjection in $\acat p$ with kernel $K$,
and let $\sh F_S \in \defs(S)$ be a deformation. Then there exists a
canonical obstruction
    \[ o(u, \sh F_S) \in ( \hh^2(\ccat, \sa, \shhom_k(\sh F_j, \sh F_i))
    \otimes_k K_{ij} ) \]
such that $o(u, \sh F_S) = 0$ if and only if there exists a
deformation $\sh F_R \in \defs(R)$ lifting $\sh F_S$ to $R$.
Moreover, if $o(u, \sh F_S) = 0$, then there is a transitive and
effective action of $( \hh^1(\ccat, \sa, \shhom_k(\sh F_j, \sh F_i))
\otimes_k K_{ij} )$ on the set of liftings of $\sh F_S$ to $R$.
\end{prop}
\begin{proof}
Let $\sh F_S \in \defs(S)$ be given. By lemma \ref{l:defm-char},
this deformation corresponds to the following data: A $k$-algebra
homomorphism $L^S(c): \sa(c) \to Q^S(c)$ for each object $c \in
\ccat$ and an element $L^S(\phi) \in Q^S(c,c')$ for each morphism
$\phi: c \to c'$ in $\ccat$, such that the conditions of lemma
\ref{l:defm-char} are satisfied. Moreover, to lift $\sh F_S$ to $R$
is the same as to lift these data to $R$.

Choose a $k$-linear section $\sigma: S \to R$ such that
$\sigma(I(S)_{ij}) \subseteq I(R)_{ij}$ and $\sigma(e_i) = e_i$ for
$1 \le i,j \le p$. Clearly, $\sigma$ induces a $k$-linear map
$Q^S(c,c') \to Q^R(c,c')$ for all $c,c'$ in $\ccat$, which we shall
denote by $Q(\sigma,c,c')$. We define $L^R(c) = Q(\sigma,c,c) \circ
L^S(c)$ for all objects $c \in \ccat$ and $L^R(\phi) =
Q(\sigma,c,c')(L^S(c,c'))$ for all morphisms $\phi: c \to c'$ in
$\ccat$. Then $L^R(c): \sa(c) \to Q^R(c)$ is a $k$-linear map which
lifts $L^S(c)$ to $R$ for all $c \in \ccat$, and we define the
obstruction
    \[ o(0,2)(c): \sa(c) \otimes_k \sa(c) \to Q^K(c) \]
by $o(0,2)(c)(a \otimes b) = L^R(c)(ab) - L^R(c)(a) \, L^R(c)(b)$
for all $a,b \in \sa(c)$. It is clear that $L^R(c)$ is a $k$-algebra
homomorphism if and only if $o(0,2)(c) = 0$. Moreover, $L^R(\phi)$
lifts $L^S(\phi)$ to $R$ for all morphisms $\phi: c \to c'$ in
$\ccat$, and we define the obstruction
    \[ o(1,1)(\phi): \sa(c) \to Q^K(c,c') \]
by $o(1,1)(\phi)(a) = L^R(\phi) \circ L^R(c)(a) -
L^R(c')(\sa(\phi)(a)) \circ L^R(\phi)$ for all $a \in \sa(c)$. It is
clear that $L^R(\phi)$ is $\sa(\phi)$-linear if and only if
$o(1,1)(\phi) = 0$. Finally, we define the obstruction
    \[ o(2,0)(\phi,\phi') \in Q^K(c,c'') \]
by $o(2,0)(\phi, \phi') = L^R(\phi') \, L^R(\phi) - L^R(\phi' \circ
\phi)$ for all morphisms $\phi: c \to c'$ and $\phi': c' \to c''$ in
$\ccat$. It is clear that $L^R$ satisfies the cocycle condition if
and only if $o(2,0) = 0$.

We see that $o = (o(0,2), o(1,1), o(2,0))$ is a $2$-cochain in
global Hochschild complex $( \hc^n(\ccat, \sa, \shhom_k(\sh F_j, \sh
F_i)) \otimes_k K_{ij} )$. A calculation shows that $o$ is a
$2$-cocycle, and that its cohomology class $o(u, \sh F_S) \in (
\hh^2(\ccat, \sa, \shhom_k(\sh F_j, \sh F_i)) \otimes_k K_{ij} )$ is
independent of the choice of $L^R(c)$ and $L^R(\phi)$. It is clear
that if there is a lifting of $\sh F_S$ to $R$, we may choose
$L^R(c)$ and $L^R(\phi)$ such that $o = 0$, hence $o(u, \sh F_S) =
0$. Conversely, assume that $o(u, \sh F_S) = 0$. Then there exists a
$1$-cochain of the form $( \epsilon, \Delta)$ with $\epsilon \in
\dc^{01}$ and $\Delta \in \dc^{10}$ such that $d(\epsilon,\Delta) =
o$. Let $L'(c) = L^R(c) + \epsilon(c)$ and $L'(\phi) = L^R(\phi) +
\Delta(\phi)$. Then $L'(c)$ is another lifting of $L^S(c)$ to $R$,
$L'(\phi)$ is another lifting of $L^S(\phi)$ to $R$, and essentially
the same calculation as above shows that the corresponding
$2$-cocycle $o' = 0$. Hence there is a lifting of $\sh F_S$ to $R$,
and this proves the first part of the proposition.

For the second part, assume that $\sh F_R$ is a lifting of $\sfam$ to
$R$. Then $\sh F_R$ is defined by liftings $L^R(c)$ and $L^R(\phi)$ to
$R$ such that the corresponding $2$-cocycle $o = 0$. Let us consider a
$1$-cochain $(\epsilon, \Delta)$, and consider the new liftings $L'(c)
= L^R(c) + \epsilon(c)$ and $L'(\phi) = L^R(\phi) + \Delta(\phi)$. From
the previous calculations, it is clear that the new $2$-cocycle $o' =
0$ if and only if $(\epsilon, \Delta)$ is a $1$-cocycle. Moreover, if
this is the case, the lifting $\sh F'_R$ defined by $L'(c)$ and
$L'(\phi)$ is equivalent to $\sh F_R$ if and only if $(\epsilon,
\Delta)$ is a $1$-coboundary, since an equivalence between $\sh F_R$
and $\sh F'_R$ must have the form $\id + \pi$ for some $0$-cochain
$\pi$ with $d(\pi) = (\epsilon, \Delta)$.
\end{proof}

We see that the obstruction $o(u, \sh F_S)$ is functorial, so it
defines an obstruction theory for the noncommutative deformation
functor $\defs: \acat p \to \sets$ by definition. If the condition
\begin{equation} \label{c:fin-coh}
    \dim_k \hh^n(\ccat, \sa, \shhom_k(\sh F_j, \sh F_i)) < \infty \text{ for }
    1 \le i,j \le p, \; n = 1,2,
\end{equation}
holds, it follows that $\defs$ has an obstruction theory with finite
dimensional cohomology $( \hh^n(\ccat, \sa, \shhom_k(\sh F_j, \sh
F_i)) )$.

\begin{thm} \label{t:defm-hull}
Let $\ccat$ be a small category, let $\sa$ be a presheaf of
$k$-algebras on $\ccat$, and let $\sfam = \{ \sh F_1, \dots, \sh F_p
\}$ be a finite family of presheaves of left $\sa$-modules on $\ccat$.
If condition (\ref{c:fin-coh}) holds, then the noncommutative
deformation functor $\defs: \acat p \to \sets$ of $\sfam$ in
$\presh(\ccat, \sa)$ has a pro-representing hull $\hull(\defs)$,
completely determined by the $k$-linear spaces $\hh^n(\ccat, \sa,
\shhom_k(\sh F_j, \sh F_i))$ for $1 \le i,j \le p$, $n = 1,2$, together
with some generalized Massey products on them.
\end{thm}

\begin{cor} \label{c:defm-iso}
Let $\ccat$ be a small category, let $\sa$ be a presheaf of
$k$-algebras on $\ccat$, and let $\sfam  = \{ \sh F_1, \dots, \sh
F_p \}$ be a finite family of presheaves of left $\sa$-modules on
$\ccat$. If (\ref{c:fin-coh}) holds, then $\hh^1(\ccat, \sa,
\shhom_k(\sh F_j, \sh F_i)) \cong t(\defs)_{ij} \cong
\ext^1_{\sa}(\sh F_j, \sh F_i)$ for $1 \le i,j \le p$.
\end{cor}

\section{Deformations of sheaves of modules} \label{s:defmsh}

Let $k$ be an algebraically closed field, and let $(X, \sa)$ be a
ringed space over $k$, i.e. a topological space $X$ together with a
sheaf $\sa$ of associative $k$-algebras on $X$. The category $\shc \sa$
of sheaves of left $\sa$-modules on $X$ is an Abelian $k$-category, and
we shall consider deformations in this category.

Let $\ccat(X)$ be the category defined in the following way: An object
in $\ccat(X)$ is an open subset $U \subseteq X$, and given objects $U,V
\in \ccat(X)$, a morphism from $U$ to $V$ in $\ccat(X)$ is an
(opposite) inclusion $U \supseteq V$. Then $\ccat(X)$ is a small
category, and we may consider $\shc \sa$ as the full subcategory of
$\presh(\ccat(X), \sa)$ consisting of exactly those presheaves of left
$\sa$-modules on $\ccat(X)$ that satisfy the sheaf axioms. It is
well-known that this subcategory is closed under direct sums, kernels
and extensions, but not under cokernels.

Let $\sfam = \{ \sh F_1, \dots, \sh F_p \}$ be a finite family of
sheaves of left $\sa$-modules on $X$. Since $\shc \sa$ is an Abelian
$k$-category and a full subcategory of $\presh(\ccat(X),\sa)$, we
may consider the noncommutative deformation functor $\defs^{sh}:
\acat p \to \sets$ of $\sfam$ as a family of sheaves, defined by
$\defs^{sh} = \defs^{\cat A}$ with $\cat A = \shc \sa$. We may also
consider the noncommutative deformation functor $\defs: \acat p \to
\sets$ of $\sfam$ as a family of presheaves, given by $\defs =
\defs^{\cat A}$ with $\cat A = \presh(\ccat(X), \sa)$, and we remark
that the natural forgetful functor $\shc \sa \to \presh(\ccat(X),
\sa)$ induces an isomorphism of noncommutative deformation functors
$\defs^{sh} \xrightarrow{\sim} \defs$. However, this observation is
not very useful for computational purposes, since the category
$\ccat(X)$ is usually too big to allow for effective computations.

Finally, we remark that $t(\defs^{sh})_{ij} \cong t(\defs)_{ij}$
since $\defs^{sh} \cong \defs$, hence $t(\defs^{sh})_{ij} \cong
\hh^1(\ccat(X),\sa,\shhom_k(\sh F_j, \sh F_i)) \cong
\ext^1_{\sa}(\sh F_j, \sh F_i)$. It does not matter if we consider
$\ext^1_{\sa}(\sh F_j, \sh F_i)$ as extensions in $\shc \sa$ or
$\presh(\ccat(X),\sa)$, since $\shc \sa$ is closed under extensions
in $\presh(\ccat(X),\sa)$.

\section{Deformations of quasi-coherent sheaves of modules}
\label{s:qcoh}

Let $k$ be an algebraically closed field, and let $(X, \sa)$ be a
ringed space over $k$. We recall that a sheaf $\sh F \in \shc \sa$ of
left $\sa$-modules on $X$ is \emph{quasi-coherent} if for every point
$x \in X$, there exists an open neighbourhood $U \subseteq X$ of $x$,
free sheaves $\sh L_0, \sh L_1$ of left $\sa \rst U$-modules on $U$,
and an exact sequence
\[ 0 \gets \sh F \rst U \gets \sh L_0 \gets \sh L_1 \]
of sheaves of left $\sa \rst U$-modules on $U$. The category $\qcoh
\sa$ of \emph{quasi-coherent sheaves} of left $\sa$-modules on $X$ is
the full subcategory of $\shc \sa$ consisting of quasi-coherent
sheaves.

The full subcategory $\qcoh \sa \subseteq \shc \sa$ is closed under
finite direct sums, but it is not clear if $\qcoh \sa$ is closed
under kernels and cokernels in general. Hence $\qcoh \sa$ is an
additive but not necessarily an Abelian $k$-category. In this
section, we give sufficient conditions for $\qcoh \sa$ to be an
exact Abelian subcategory of $\shc \sa$, and consider deformations
in the category $\qcoh \sa$ in these cases.

Let us consider the global sections functor $\Gm(X,-): \shc \sa \to
\mcat A$, where we write $A = \Gm(X,\sa)$. This functor is left exact,
and we denote its right derived functors by $\ch^{\cpd}(X,-) = R^{\cpd}
\Gm(X,-)$. We say that $X$ is \emph{$\sa$-affine} if the following
conditions hold:
\begin{enumerate}
\item \label{c:aff-i}
$\Gm(X,-)$ induces an equivalence of categories $\qcoh \sa \to \mcat
A$,
\item \label{c:aff-ii}
$\ch^n(X,\sh F) = 0$ for all $\sh F \in \qcoh \sa$ and for all
integers $n \ge 1$.
\end{enumerate}
Moreover, we say that an open subset $U \subseteq X$ is $\sa$-affine if
$U$ is $\sa \rst U$-affine, and that an open cover $\oc U$ of $X$ is
$\sa$-affine if $U$ is $\sa$-affine for any $U \in \oc U$.

We say that a full subcategory $\cat C$ of an Abelian $k$-category
$\cat A$ is \emph{thick} if $\cat C \subseteq \cat A$ is an exact
Abelian subcategory that is closed under extensions, and recall that
\begin{enumerate}
\item If $X_1 \to X_2 \to Y \to X_3 \to X_4$ is an exact sequence in
$\cat A$ with $X_i \in \cat C$ for $1 \le i \le 4$, then $Y \in \cat
C$.
\end{enumerate}
is a necessary and sufficient condition for $\cat C \subseteq \cat A$
to be a thick subcategory.

\begin{prop} \label{p:abcat}
If $X$ has an $\sa$-affine open cover, then $\qcoh \sa \subseteq \shc
\sa$ is a thick subcategory.
\end{prop}
\begin{proof}
Using the condition (1) above, we see that $\qcoh \sa$ is a thick
subcategory of $\shc \sa$ if and only this holds locally. We may
therefore assume that $X$ is $\sa$-affine. This implies that $\qcoh
\sa$ has kernels and cokernels, since this holds for $\mcat A$. We must
show that $\qcoh \sa$ is closed under extensions. Let $0 \to \sh F \to
\sh G \to \sh H \to 0$ be an exact sequence in $\shc \sa$ with $\sh F,
\sh H$ in $\qcoh \sa$. Then $\ch^1(X, \sh F) = 0$, and $\Gm(X,-): \shc
\sa \to \mcat A$ is right adjoint to $\qc X: \mcat A \to \qcoh \sa$,
the quasi-inverse of $\Gm(X,-)$ restricted to $\qcoh \sa$. So there is
a commutative diagram
\[
    \xymatrix{
    0 \ar[r] & \qc X \Gm(X, \sh F) \ar[r] \ar[d] & \qc X
    \Gm(X, \sh G) \ar[r] \ar[d] & \qc X \Gm(X, \sh H) \ar[r]
    \ar[d] & 0 \\
    0 \ar[r] & \sh F \ar[r] & \sh G \ar[r] & \sh H \ar[r] & 0
    }
\]
in $\shc \sa$ with exact rows. The left and right vertical arrows
are isomorphisms, so the middle vertical arrow is an isomorphism as
well.
\end{proof}

Let $\oc U$ be an $\sa$-affine open cover of $X$, viewed as a small
subcategory of $\ccat(X)$, and let $\pi: \qcoh \sa \to \presh(\oc U,
\sa)$ be the natural forgetful functor. For any finite family $\sfam
= \{ \sh F_1, \dots, \sh F_p \}$ of quasi-coherent sheaves of left
$\sa$-modules on $X$, we may consider the noncommutative deformation
functor $\defs^{qc}: \acat p \to \sets$ of $\sfam$ as a family of
quasi-coherent sheaves, defined by $\defs^{qc} = \defs^{\cat A}$
with $\cat A = \qcoh \sa$. We may also consider the noncommutative
deformation functor $\defs^{\oc U}: \acat p \to \sets$ of $\sfam$ as
a family of presheaves on $\oc U$, defined by $\defs^{\oc U} =
\defs^{\cat A}$ with $\cat A = \presh(\oc U, \sa)$. We remark that
$\pi$ induces a morphism $\defs^{qc} \to \defs^{\oc U}$ of
noncommutative deformation functors, but not necessarily an
isomorphism.

An open cover $\oc U$ of $X$ is \emph{good} if any finite intersection
$V = U_1 \cap U_2 \cap \dots \cap U_r$ with $U_i \in \oc U$ for $1 \le
i \le r$ can be covered by open subsets $W \subseteq V$ with $W \in \oc
U$. In particular, any open cover closed under finite intersections is
good.

\begin{prop} \label{p:defm-eqv}
If $\oc U$ is a good $\sa$-affine open cover of $X$, then the
forgetful functor $\pi: \qcoh \sa \to \presh(\oc U, \sa)$ induces an
isomorphism $\defs^{qc} \to \defs^{\oc U}$ of deformation functors
for any finite family $\sfam$ of quasi-coherent left $\sa$-modules
on $X$.
\end{prop}
\begin{proof}
Clearly, $\pi$ induces a morphism of noncommutative deformation
functors, and it is enough to show that the induced map of sets
$\pi^*_R: \defs^{qc}(R) \to \defs^{\oc U}(R)$ is a bijection for any
$R \in \acat p$. If $X$ is $\sa$-affine and $\oc U = \{ X \}$, then
$\presh(\oc U, \sa)$ is naturally equivalent to $\mcat A$, so $\pi:
\qcoh \sa \to \presh(\oc U, \sa)$ is an equivalence of categories,
and this implies that $\pi_R$ is a bijection for any $R \in \acat
p$. In the general case, let $\sh F_R \in \defs^{\oc U}(R)$. Then
$\sh F_R(U)$ is a deformation of the family $\{ \sh F_1(U), \dots,
\sh F_p(U) \}$ in $\mcat{\sa(U)}$ to $R$ for any $U \in \oc U$. By
the result in the $\sa$-affine case, we can find a deformation $\sh
F_R^U$ of the family $\{ \sh F_1 \rst U, \dots, \sh F_p \rst U \}$
in $\qcoh{\sa \rst U}$ to $R$ that is compatible with $\sh F_R(U)$.
We remark that if $V \subseteq U$ is an inclusion in $\oc U$, then
there is a natural isomorphism $\sh F_R^U \rst V \to \sh F_R^V$ of
sheaves of left $\sa$-modules on $V$, since $\sfam$ is a family of
quasi-coherent sheaves of $\sa$-modules on $X$. We must glue the
local deformations $\sh F_R^U$ to a deformation $\sh F_R$ of the
family $\sfam$ to $R$ in $\qcoh \sa$, and this is clearly possible
since $\oc U$ is a good open cover of $X$.
\end{proof}

\begin{thm} \label{t:defmqc-hull}
Let $(X,\sa)$ be a ringed space over $k$, let $\oc U$ be a good
$\sa$-affine open cover of $X$, and let $\sfam = \{ \sh F_1, \dots, \sh
F_p \}$ be a finite family of quasi-coherent left $\sa$-modules on $X$.
If $\dim_k \hh^n(\oc U, \sa, \shhom_k(\sh F_j, \sh F_i)) < \infty$ for
$1 \le i,j \le p, \; n = 1,2$, then the noncommutative deformation
functor $\defs^{qc}: \acat p \to \sets$ of $\sfam$ in $\qcoh \sa$ has a
pro-representing hull $\hull(\defs^{qc})$, completely determined by the
$k$-linear spaces $\hh^n(\oc U, \sa, \shhom_k(\sh F_j, \sh F_i))$ for
$1 \le i,j \le p$, $n = 1,2$, together with some generalized Massey
products on them.
\end{thm}

We remark that if the conditions of theorem \ref{t:defmqc-hull}
holds, then $t(\defs^{qc})_{ij} \cong t(\defs^{\oc U})_{ij}$ since
$\defs^{qc} \cong \defs^{\oc U}$, hence $t(\defs^{qc})_{ij} \cong
\hh^1(\oc U,\sa,\shhom_k(\sh F_j, \sh F_i)) \cong \ext^1_{\sa}(\sh
F_j, \sh F_i)$. It does not matter if we consider $\ext^1_{\sa}(\sh
F_j, \sh F_i)$ as extensions in $\presh(\oc U,\sa)$ or in $\qcoh
\sa$, since the forgetful functor $\qcoh \sa \to \presh(\oc U, \sa)$
preserves extensions by the argument in the proof of proposition
\ref{p:defm-eqv}.

\section{Quasi-coherent ringed schemes} \label{s:imp-ex}

Let $k$ be an algebraically closed field. We shall consider some
important examples of ringed spaces $(X, \sa)$ over $k$ such that $X$
has a good $\sa$-affine open cover $\oc U$. It is often possible to
choose $\oc U$ to be a finite cover, and this is important for
effective computations of noncommutative deformations using presheaf
methods.

\begin{exc} \label{e:sch}
Let $(X, \reg X)$ be a scheme over $k$. If $U \subseteq X$ is an
open affine subscheme of $X$, then $U$ is $\reg X$-affine by
Hartshorne \cite{hart77}, corollary II.5.5 and Grothendieck
\cite{grot61ii}, theorem 1.3.1. Hence any open affine cover of $X$
is an $\reg X$-affine open cover. If $X$ is separated over $k$, then
any finite intersection of open affine subschemes of $X$ is affine.
Hence if $X$ is quasi-compact and separated over $k$, then there is
a finite $\reg X$-affine open cover of $X$ closed under
intersections.
\end{exc}

\begin{exc} \label{e:ass}
A \emph{ringed scheme} over $k$ is a ringed space $(X, \sa)$ over
$k$ defined by a scheme $(X, \reg X)$ over $k$ and a morphism $i:
\reg X \to \sa$ of sheaves of associative $k$-algebras on $X$. A
\emph{quasi-coherent ringed scheme} over $k$ is a ringed scheme $(X,
\sa)$ over $k$ such that $\sa$ is quasi-coherent as a left and right
$\reg X$-module. The notion of quasi-coherent ringed schemes was
introduced in Yekutieli, Zhang \cite{yek-zha05}, and the following
result follows from Yekutieli, Zhang \cite{yek-zha05}, corollary
5.13:

\begin{lem} \label{l:loc}
A ringed scheme $(X, \sa)$ over $k$ is quasi-coherent if and only if
the morphism $\sa(U) \to \sa(D(f))$ is a ring of fractions with
respect to $S = \{ f^n: n \ge 0 \}$ for any open affine subscheme $U
\subseteq X$ and any $f \in \reg X(U)$.
\end{lem}

For any quasi-coherent ringed scheme $(X, \sa)$ over $k$, a left
$\sa$-module $\sh F$ is quasi-coherent if and only if it is
quasi-coherent as a left $\reg X$-module. This follows from
Grothendieck \cite{grot60}, proposition 9.6.1 when $\sa$ is a sheaf
of commutative rings on $X$, and the proof can easily be extended to
the noncommutative case.

\begin{lem}
Let $(X, \sa)$ be a quasi-coherent ringed scheme over $k$. If $U
\subseteq X$ is an open affine subscheme of $X$, then $U$ is
$\sa$-affine.
\end{lem}
\begin{proof}
Write $A = \sa(U)$, and consider $\Gm(U,-): \qcoh{\sa \rst U} \to \mcat
A$. We claim that $\Gm(U,-)$ is an equivalence of categories. By the
comments preceding the lemma, $\qcoh{\sa \rst U}$ can be considered as
a subcategory of $\qcoh{\reg U}$, and $\Gm(U,-): \qcoh{\reg U} \to
\mcat O$ is an equivalence of categories with $O = \reg X(U)$. So the
claim follows from lemma \ref{l:loc}. Finally, $\ch^n(U, \sh F) = 0$
for any integer $n \ge 1$ and any $\sh F \in \qcoh{\sa \rst U}$ by the
above comments.
\end{proof}

Let $(X, \sa)$ be a quasi-coherent ringed scheme over $k$. Then any
open affine cover of $X$ is an $\sa$-affine open cover of $X$. If $X$
is quasi-compact and separated over $k$, then there is a finite
$\sa$-affine open cover of $X$ closed under intersections.
\end{exc}

\begin{exc} \label{e:bbdiff}
Let $(X, \sa)$ be a quasi-coherent ringed scheme over $k$, and
assume that $\car(k) = 0$. We say that $\sa$ is a \emph{D-algebra},
and that $(X, \sa)$ is a \emph{D-scheme}, if the following condition
holds: For any open subset $U \subseteq X$ and for any section $a
\in \sa(U)$, there exists an integer $n \ge 0$ (depending on $a$)
such that
\[ [ \dots [ \; [a, f_1], f_2] \dots , f_n] = 0 \]
for all sections $f_1, \dots, f_n \in \reg X(U)$, where $[a,f] = a f
- f a$ is the usual commutator for all $a \in \sa(U), \; f \in \reg
X(U)$. The notion of D-schemes was considered in Beilinson,
Bernstein \cite{bei-ber93}, and most quasi-coherent ringed schemes
that appear naturally are D-schemes. We give some important examples
of D-schemes below.
\end{exc}

\begin{exc} \label{e:diff}
Let $(X, \reg X)$ be a scheme over $k$, and assume that $\car(k) =
0$. For any sheaf $\sh F$ of $\reg X$-modules, we denote by $\diff
F$ the sheaf of $k$-linear \emph{differential operators} on $\sh F$,
see Grothendieck \cite{grot67}, section 16.8. By definition, $\diff
F$ is a sheaf of associative $k$-algebras on $X$, equipped with a
morphism $i: \reg X \to \diff F$ of sheaves of rings, and $\diff F$
is clearly a D-algebra on $X$ if and only if $\diff F$ is
quasi-coherent as a left and right $\reg X$-module. By Beilinson,
Bernstein \cite{bei-ber93}, example 1.1.6, this is the case if $\sh
F$ is a coherent $\reg X$-module.

Let $X$ be locally Noetherian, and consider the sheaf $\sh D_X =
\diff{\reg X}$ of $k$-linear \emph{differential operators} on $X$.
Since $\reg X$ is a coherent sheaf of rings, it follows that $\sh D_X$
is a D-algebra on $X$.

We remark that there are some examples of schemes over $k$ that are
$\sh D_X$-affine but not affine. For instance, this holds for the
projective space $X = \pp^n$ for all integers $n \ge 1$, see
Beilinson, Bernstein \cite{bei-ber81}. It also holds for the
weighted projective space $X = \pp(a_1, \dots, a_n)$, see Van den
Bergh \cite{vdbe91}.
\end{exc}

\begin{exc} \label{e:liealg}
Let $(X, \reg X)$ be a separated scheme of finite type over $k$, and
assume that $\car(k) = 0$. A \emph{Lie algebroid} of $X$ is a
quasi-coherent $\reg X$-module $\g$ with a $k$-Lie algebra structure,
together with a morphism $\tau: \g \to \shder_k(\reg X)$ of sheaves of
$\reg X$-modules and of $k$-Lie algebras, such that
    \[ [g, f \cdot h] = f [g, h] + \tau_g(f) \cdot h \]
for any open subset $U \subseteq X$ and any sections $g, h \in
\g(U), \; f \in \reg X(U)$. The notion of Lie algebroids in
algebraic geometry was considered in Beilinson, Bernstein
\cite{bei-ber93}.

For any sheaf $\sh F$ of $\reg X$-modules, an integrable
\emph{$\g$-connection} on $\sh F$ is a morphism $\nabla: \g \to
\shend_k(\sh F)$ of sheaves of $\reg X$-modules and of $k$-Lie
algebras such that
\[ \nabla_U(g)(f m) = f \; \nabla_U(g)(m) + g(f) \; m \]
for any open subset $U \subseteq X$ and any sections $f \in \reg
X(U), \; g \in \g(U), \; m \in \sh F(U)$. The quasi-coherent sheaves
of $\reg X$-modules with integrable $\g$-connections form an Abelian
$k$-category, and there is a universal enveloping D-algebra
$\ue(\g)$ of $\g$ such that this category is equivalent to
$\qcoh{\ue(\g)}$, see Beilinson, Bernstein \cite{bei-ber93}. In
particular, $(X, \ue(\g))$ is a D-scheme.

The tangent sheaf $\theta_X = \shder_k(\reg X)$ of $X$ is a Lie
algebroid of $X$ in a natural way, and $\ue(\theta_X)$ is a subsheaf
of the sheaf $\sh D_X$ of $k$-linear differential operators on $X$.
If $X$ is a smooth irreducible quasi-projective variety over $k$,
then $\ue(\theta_X) = \sh D_X$.
\end{exc}

\section{Calculations for D-modules on elliptic curves}
\label{s:example}

Let $k$ be an algebraically closed field of characteristic $0$, and
let $X$ be a smooth irreducible variety over $k$ of dimension $d$.
Then the sheaf $\sh D_X$ of $k$-linear differential operators on $X$
is a D-algebra on $X$. We consider the noncommutative deformations
of $\reg X$ as a quasi-coherent left $\sh D_X$-module via the
natural left action of $\sh D_X$ on $\reg X$. As an example, we
compute the pro-representing hull of $\reg X$ as a quasi-coherent
left $\sh D_X$-module when $X$ is an elliptic curve, see also
Eriksen \cite{erik07}.

Let $\oc U$ be an open affine cover of $X$. Then $U \subseteq X$ is
a smooth, irreducible affine variety over $k$ of dimension $d$ for
all $U \in \oc U$. It is well-known that $\sh D_X(U)$ is a simple
Noetherian ring of global dimension $d$ and that $\reg X(U)$ is a
simple left $\sh D_X(U)$-module, see Smith, Stafford
\cite{smi-sta88}. Hence $\shext^q_{\sh D_X}(\reg X, \reg X) = 0$ for
$q \ge d+1$ and $\shend_{\sh D_X}(\reg X) = \underline k$. If $X$ is
curve, then the spectral sequence in proposition \ref{p:glh-ss}
degenerates, and
\begin{align*}
    \hh^n(\oc U, \sh D_X, \shend_k(\reg X)) & \cong \ch^{n-1}(\oc U,
    \shext^1_{\sh D_X}(\reg X, \reg X)) \text{ for } n \ge 1 \\
    \hh^0(\oc U, \sh D_X, \shend_k(\reg X)) & \cong k
\end{align*}

Let $X \subseteq \pp^2$ be the irreducible projective plane curve
given by the homogeneous equation $f = 0$, where $f = y^2 z - x^3 -
a x z^2 - b z^3$ for fixed parameters $(a,b) \in k^2$. We assume
that $\Delta = 4a^3 + 27b^2 \neq 0$, so that $X$ is smooth and
therefore an elliptic curve over $k$. We choose an open affine cover
$\oc U = \{ U_1, U_2, U_3 \}$ of $X$ closed under intersections,
given by $U_1 = D_+(y)$, $U_2 = D_+(z)$ and $U_3 = U_1 \cap U_2$. We
shall compute $\shext^1_{\sh D_X}(\reg X, \reg X)$ and
$\ch^{n-1}(\oc U, \shext^1_{\sh D_X}(\reg X, \reg X))$ for $n =
1,2$.

Let $A_i = \reg X(U_i)$ and $D_i = \sh D_X(U_i)$ for $i=1,2,3$. We
see that $A_1 \cong k[x,z]/(f_1)$ and $A_2 \cong k[x,y]/(f_2)$,
where $f_1 = z - x^3 - axz^2 - bz^3$ and $f_2 = y^2 - x^3 - ax - b$.
Moreover, we have that $\der_k(A_i) = A_i \partial_i$ and $D_i = A_i
\langle \partial_i \rangle$ for $i=1,2$, where
\begin{align*}
    \partial_1 & = (1 - 2axz - 3bz^2) \; \partial / \partial x +
    (3x^2 + az^2) \; \partial / \partial z \\
    \partial_2 & = -2y \; \partial / \partial x - (3x^2 + a) \;
    \partial / \partial y
\end{align*}
On the intersection $U_3 = U_1 \cap U_2$, we choose an isomorphism
$A_3 \cong k[x,y,y^{-1}]/(f_3)$ with $f_3 = f_2$, and see that
$\der_k(A_3) = A_3 \partial_3$ and $D_3 = A_3 \langle \partial_3
\rangle$ for $\partial_3 =
\partial_2$. The restriction maps of $\reg X$ and $\sh D_X$, considered
as presheaves on $\oc U$, are given by
    \[ x \mapsto xy^{-1}, \; z \mapsto y^{-1}, \; \partial_1 \mapsto
    \partial_2 \]
for the inclusion $U_1 \supseteq U_3$, and the natural localization map
for $U_2 \supseteq U_3$. Finally, we find a free resolution of $A_i$ as
a left $D_i$-module for $i=1,2,3$, given by
\begin{align*}
0 \gets A_i & \gets D_i \xleftarrow{\cdot \partial_i} D_i \gets 0
\end{align*}
and use this to compute $\ext^1_{D_i}(A_i, A_j) \cong
\coker(\partial_i \rst{U_j}: A_j \to A_j)$ for all $U_i \supseteq
U_j$ in $\oc U$. We see that $\ext^1_{D_i}(A_i, A_3) \cong
\coker(\partial_3: A_3 \to A_3)$ is independent of $i$, and find the
following $k$-linear bases for $\ext^1_{D_i}(A_i,A_j)$:

\vspace{5mm}
\begin{center}
\begin{tabular}{|l|l|l|}
    \hline
    & $a \neq 0:$ & $a = 0:$ \\
    \hline
    $U_1 \supseteq U_1$ & $1, z, z^2, z^3$ & $1, z, x, xz$ \\
    $U_2 \supseteq U_2$ & $1, y^2$ & $1, x$ \\
    $U_3 \supseteq U_3$ & $x^2 y^{-1}, 1, y^{-1}, y^{-2}, y^{-3}$
    & $x^2 y^{-1}, 1, y^{-1}, x, x y^{-1}$ \\
    \hline
\end{tabular}
\end{center}
\vspace{5mm}

\noindent The functor $\shext^1_{\sh D_X}(\reg X, \reg X):
\morcat{\oc U} \to \mcat k$ defines the following diagram in $\mcat
k$, where the maps are induced by the restriction maps on $\reg X$:
\[
\xymatrix{ \ext^1_{D_1}(A_1,A_1) \ar[d] & & \ext^1_{D_2}(A_2,A_2)
\ar[d] \\
\ext^1_{D_1}(A_1,A_3) & \ext^1_{D_3}(A_3,A_3) \ar@{=}[l] \ar@{=}[r] &
\ext^1_{D_2}(A_2,A_3) }
\]
We use that $15 y^2 = \Delta \; y^{-2}$ in $\ext^1_{D_3}(A_3,A_3)$
when $a \neq 0$ and that $-3b \; xy^{-2} = x$ in
$\ext^1_{D_3}(A_3,A_3)$ when $a = 0$ to describe these maps in the
given bases, and compute $\ch^{n-1}(\oc U, \ext^1_{\sh D_X}(\reg X,
\reg X))$ for $n = 1,2$ using the resolving complex $\dc^{\cpd}(\oc
U, -)$. We find the following $k$-linear bases:

\vspace{5mm}
\begin{center}
\begin{tabular}{|l|l|l|}
    \hline
    & $a \neq 0:$ & $a = 0:$ \\
    \hline
    $n=1$ & $\xi_1 = (1,1,1), \; \xi_2 = (\Delta z^2, 15 y^2, \Delta
    y^{-2})$ & $\xi_1 = (1,1,1), \; \xi_2 = (-3b \; xz, x, x)$ \\
    $n=2$ & $\omega = (0,0,0,0,6a x^2 y^{-1})$ & $\omega =
    (0,0,0,0,x^2 y^{-1})$ \\
    \hline
\end{tabular}
\end{center}
\vspace{5mm}

\noindent We recall that $\xi_1, \xi_2$ and $\omega$ are represented
by cocycles of degree $p=0$ and $p=1$ in the resolving complex
$\dc^{\cpd}(\oc U, \shext^1_{\sh D_X}(\reg X, \reg X))$, where
    \[ \dc^p(\oc U, \shext^1_{\sh D_X}(\reg X, \reg X)) = \prod_{U_0
    \supseteq \dots \supseteq U_p} \shext^1_{\sh D_X}(\reg X,\reg X)(U_0
    \supseteq U_p) \]
and the product is indexed by $\{ U_1 \supseteq U_1, U_2 \supseteq
U_2, U_3 \supseteq U_3 \}$ when $p = 0$, and $\{ U_1 \supseteq U_1,
U_2 \supseteq U_2, U_3 \supseteq U_3, U_1 \supseteq U_3, U_2
\supseteq U_3 \}$ when $p = 1$.

This proves that the noncommutative deformation functor $\defm{\reg
X}: \acat 1 \to \sets$ of the left $\sh D_X$-module $\reg X$ has
tangent space $\hh^1(\oc U, \sh D_X, \shend_k(\reg X)) \cong k^2$
and obstruction space $\hh^2(\oc U, \sh D_X, \shend_k(\reg X)) \cong
k$ for any elliptic curve $X$ over $k$, and a pro-representing hull
$H = k {\ll} t_1,t_2 {\gg} /(F)$ for some noncommutative power
series $F \in k {\ll} t_1,t_2 {\gg}$.

We shall compute the noncommutative power series $F$ and the versal
family $\sh F_H \in \defm{\reg X}(H)$ using the obstruction
calculus. We choose base vectors $t_1^*, t_2^*$ in $\hh^1(\oc U,
\sa, \shend_k(\reg X))$, and representatives $(\psi_l, \tau_l) \in
\dc^{01} \oplus \dc^{10}$ of $t_l^*$ for $l = 1,2$, where $\dc^{pq}
= \dc^p(\oc U, \hc^q(\sa, \shend_k(\reg X)))$. We may choose
$\psi_l(U_i)$ to be the derivation defined by
\begin{equation*}
\psi_l(U_i)(P_i) = \begin{cases} 0 & \text{if $P_i \in A_i$} \\
\xi_l(U_i) \cdot \id_{A_i} & \text{if $P_i = \partial_i$}
\end{cases}
\end{equation*}
for $l = 1,2$ and $i = 1,2,3$, and $\tau_l(U_i \supseteq U_j)$ to be
the multiplication operator in $\hmm_{A_i}(A_i,A_j) \cong A_j$ given
by $\tau_1 = 0$, $\tau_2(U_i \supseteq U_i) = 0$ for $i = 1,2,3$ and

\vspace{5mm}
\begin{center}
\begin{tabular}{|l|l|l|}
    \hline
    $a \neq 0:$ & $a = 0:$ \\
    \hline
    $\tau_2(U_1 \supseteq U_3) = 0$ & $\tau_2(U_1 \supseteq U_3) =
    x^2y^{-1}$ \\
    $\tau_2(U_2 \supseteq U_3) = -4a^2y^{-1} - 3xy + 9bxy^{-1} -
    6ax^2y^{-1}$ & $\tau_2(U_2 \supseteq U_3) = 0$ \\
    \hline
\end{tabular}
\end{center}
\vspace{5mm}

The restriction of $\defm{\reg X}: \acat 1 \to \sets$ to $\acat
1(2)$ is represented by $(H_2, \sh F_{H_2})$, where $H_2 = k \langle
t_1,t_2 \rangle /(t_1,t_2)^2$ and the deformation $\sh F_{H_2} \in
\defm{\reg X}(H_2)$ is defined by $\sh F_{H_2}(U_i) = A_i \otimes_k
H_2$ as a right $H_2$-module for $i = 1,2,3$, with left $D_i$-module
structure given by
\begin{equation*}
P_i (m_i \otimes 1) = P_i(m_i) \otimes 1 + \psi_1(U_i)(P_i)(m_i)
\otimes t_1 + \psi_2(U_i)(P_i)(m_i) \otimes t_2
\end{equation*}
for $i = 1,2,3$ and for all $P_i \in D_i, \; m_i \in A_i$, and with
restriction map for the inclusion $U_i \supseteq U_j$ given by
\begin{equation*}
m_i \otimes 1 \mapsto m_i \rst{U_j} \otimes 1 + \tau_2(U_i \supseteq
U_j) \; m_i \rst{U_j} \otimes t_2
\end{equation*}
for $i=1,2, \; j = 3$ and for all $m_i \in A_i$.

Let us attempt to lift the family $\sh F_{H_2} \in \defm{\reg
X}(H_2)$ to $R = k{\ll}t_1,t_2{\gg}/(t_1,t_2)^3$. We let $\sh
F_R(U_i) = A_i \otimes_k R$ as a right $R$-module for $i = 1,2,3$,
with left $D_i$-module structure given by
\begin{equation*}
P_i (m_i \otimes 1) = P_i(m_i) \otimes 1 + \psi_1(U_i)(P_i)(m_i)
\otimes t_1 + \psi_2(U_i)(P_i)(m_i) \otimes t_2
\end{equation*}
for $i = 1,2,3$ and for all $P_i \in D_i, \; m_i \in A_i$, and with
restriction map for the inclusion $U_i \supseteq U_j$ given by
\begin{equation*}
m_i \otimes 1 \mapsto m_i \rst{U_j} \otimes 1 + \tau_2(U_i \supseteq
U_j) \; m_i \rst{U_j} \otimes t_2 + \frac{\tau_2(U_i \supseteq
U_j)^2}{2} \; m_i \rst{U_j} \otimes t_2^2
\end{equation*}
for $i=1,2, \; j = 3$ and for all $m_i \in A_i$. We see that $\sh
F_R(U_i)$ is a left $\sh D_X(U_i)$-module for $i = 1,2,3$, and that
$t_1 t_2 - t_2 t_1 = 0$ is a necessary and sufficient condition for
$\sh D_X$-linearity of the restriction maps for the inclusions $U_1
\supseteq U_3$ and $U_2 \supseteq U_3$. This implies that $\sh F_R$
is not a lifting of $\sh F_{H_2}$ to $R$. But if we consider the
quotient $H_3 = R/(t_1 t_2 - t_2 t_1)$, we see that the family $\sh
F_{H_3} \in \defm{\reg X}(H_3)$ induced by $\sh F_R$ is a lifting of
$\sh F_{H_2}$ to $H_3$.

In fact, we claim that the restriction of $\defm{\reg X}: \acat 1
\to \sets$ to $\acat 1(3)$ is represented by $(H_3, \sh F_{H_3})$.
One way to prove this is to show that it is not possible to find any
lifting $\sh F'_R \in \defm{\reg X}(R)$ of $\sh F_{H_2}$ to $R$.
Another approach is to calculate the cup products $< t^*_i, t^*_j>$
in global Hochschild cohomology for $i,j = 1,2$, and this gives
\begin{align*}
    < t^*_1, t^*_2> = o^*, \; & < t^*_2, t^*_1> = - o^* &
    \text{ for } a \neq 0 \\
    < t^*_1, t^*_2> = o^*, \; & < t^*_2, t^*_1> = - o^* &
    \text{ for } a = 0
\end{align*}
where $o^* \in \hh^2(\oc U, \sh D_X, \shend_k(\reg X))$ is the base
vector corresponding to $\omega$. Since all other cup products
vanish, this implies that $F = t_1 t_2 - t_2 t_1 \mod (t_1,t_2)^3$.

Let $H = k{\ll}t_1,t_2{\gg}/(t_1t_2 - t_2 t_1)$. We shall show that it
is possible to find a lifting $\sh F_H \in \defm{\reg X}(H)$ of $\sh
F_{H_3}$ to $H$. We let $\sh F_\hull(U_i) = A_i \widehat{\otimes}_k H$
as a right $H$-module for $i = 1,2,3$, with left $D_i$-module structure
given by
\begin{equation*}
P_i (m_i \otimes 1) = P_i(m_i) \otimes 1 + \psi_1(U_i)(P_i)(m_i)
\otimes t_1 + \psi_2(U_i)(P_i)(m_i) \otimes t_2
\end{equation*}
for $i = 1,2,3$ and for all $P_i \in \sh D_i, \; m_i \in A_i$, and
with restriction map for the inclusion $U_i \supseteq U_j$ given by
\begin{equation*}
    m_i \otimes 1 \mapsto \sum_{n=0}^{\infty} \frac{\tau_2(U_i
    \supseteq U_j)^n}{n!} \; m_1 \rst{U_j} \otimes t_2^n
    = \exp(\tau_2(U_i \supseteq U_j) \otimes t_2) \cdot (m_1
    \rst{U_j} \otimes 1)
\end{equation*}
for $i=1,2, \; j = 3$ and for all $m_i \in A_i$. This implies that
$(H, \sh F_H)$ is the pro-representing hull of $\defm{\reg X}$, and
that $F = t_1 t_2 - t_2 t_1$. We remark that the versal family $\sh
F_H$ does not admit an algebraization, i.e. an algebra
$H_{\text{alg}}$ of finite type over $k$ such that $H$ is a
completion of $H_{\text{alg}}$, together with a deformation in
$\defm{\reg X}(H_{\text{alg}})$ that induces the versal family $\sh
F_H \in \defm{\reg X}(H)$.

\bibliographystyle{amsplain}
\bibliography{main}

\end{document}